\documentclass[12pt,twoside]{article}

\usepackage{xcolor}
\usepackage[english]{babel}
\usepackage{blindtext}
\usepackage[utf8]{inputenc}
\usepackage[T1]{fontenc}
\usepackage{graphicx} 
\usepackage{amssymb}
\usepackage{amsmath}
\usepackage{subfigure}
\usepackage{amsthm}
\usepackage{relsize}
\usepackage{blindtext}
\usepackage{titlesec}
\usepackage{comment}
\usepackage[toc,page]{appendix}
\usepackage[a4paper, width=150mm, top=25mm, bottom=25mm, bindingoffset=6mm]{geometry}
\usepackage{mathtools}
\usepackage[hidelinks]{hyperref}
\usepackage{authblk}
\usepackage{enumerate}

\usepackage[backend=biber,backref, giveninits=true,isbn=false, url=false, maxbibnames=50]{biblatex}

\usepackage{fancyhdr}

\usepackage[nameinlink]{cleveref}

\crefname{lstlisting}{listing}{listings}
\Crefname{lstlisting}{Listing}{Listings}

\usepackage{hyperref}
\hypersetup{
	pdfauthor={von der Mosel and Strzelecki},
	pdftitle={TP-minimizing spheres},
	colorlinks=true,
	linktoc=all,  
	linkcolor=blue,
	urlcolor=blue
}

\newtheorem{thm}{Theorem}[section]
\newtheorem{defn}[thm]{Definition}
\newtheorem{rmk}[thm]{Remark}
\newtheorem{prop}[thm]{Proposition}
\newtheorem{lem}[thm]{Lemma}
\newtheorem{nota}[thm]{Notation}
\newtheorem{cor}[thm]{Corollary}

\numberwithin{equation}{section}

\definecolor{brightviolet}{HTML}{8A2BE2}

\newcommand{\omitted}[1]{}

\newcommand{\R}{\mathbb{R}}

\DeclareFontFamily{U}{matha}{}
\DeclareFontShape{U}{matha}{m}{n}{
  <5><6><7><8><9><10>gen*matha
  <10.95><12><14.4><17.28><20.74><24.88>matha10
}{}
\DeclareSymbolFont{matha}{U}{matha}{m}{n}

\makeatletter
\newcommand{\preprintContacts}{}

\providecommand{\address}[2][]{}
\renewcommand{\address}[2][]{%
	\g@addto@macro\preprintContacts{%
		\par\noindent(#1)\hspace{0.2em}\textsc{#2}\par
	}%
}

\providecommand{\email}[2][]{}
\renewcommand{\email}[2][]{%
	\g@addto@macro\preprintContacts{%
		\noindent \textit{Email address}: \texttt{"\detokenize{#2}"}\par\medskip
	}%
}
\makeatother

\AtEndDocument{%
	\par\bigskip
	{\fontsize{9pt}{9pt}\selectfont\preprintContacts\par}%
}

\title{Round hyperspheres uniquely minimize most tangent-point energies}

\author{Alexander Dohmen} 
\address[Alexander Dohmen]{
	RWTH Aachen University,
	Institute for Mathematics,
	Templergraben~55,
	52062~Aachen, Germany}
\email[Alexander Dohmen]{dohmen@instmath.rwth-aachen.de}

\author{Heiko von der Mosel}
\address[Heiko von der Mosel]{
	RWTH Aachen University,
	Institute for Mathematics,
	Templergraben~55,
	52062~Aachen, Germany}
\email[Heiko von der Mosel]{heiko@instmath.rwth-aachen.de}

\affil{Institute for Mathematics, RWTH Aachen University}

\date{\today}

\begin{document}
\def\TP{{\rm{TP}}}

\begingroup
\let\newpage\relax
\let\clearpage\relax
\let\cleardoublepage\relax

\maketitle

\begin{abstract}
    We show that round hyperspheres uniquely minimize a broad class of
    generalized tangent-point energies among all boundaries of bounded Caccioppoli sets sharing the same perimeter or enclosed volume. This
    includes closed $C^1$-hypersurfaces. Equivalently, we establish sharp lower bounds for the associated scale-invariant ratios of energy to surface measure and energy to volume, respectively, with equality precisely for round spheres.  
\end{abstract}

\clearpage

\endgroup

\section{Introduction}

We consider the tangent-point energies $\TP_q$ and generalized tangent-point energies $\TP^{(p,q)}$ on the domain of boundaries $\Sigma$ of bounded \emph{Caccioppoli sets} $\Omega \subset \mathbb{R}^{d+1}$. These sets are defined as sets of finite perimeter $\textnormal{Per}(\Omega)$, which is the total variation of the characteristic function of $\Omega$ in $\R^{d+1}$. Caccioppoli sets have a \textit{measure-theoretic outward unit-normal vector field}  $n=n_\Omega: \Sigma ^\ast \rightarrow \mathbb{R}^{d+1}$, defined on the $d-$rectifiable \emph{reduced boundary} $\Sigma^\ast :=\partial ^\ast \Omega \subseteq \Sigma$ of $\Omega$ satisfying $\textnormal{Per}(\Omega)= \mathcal{H}^{d}(\Sigma ^\ast)$. Here, $\mathcal{H}^d$ denotes the $d$-dimensional Hausdorff measure. For example, every 
open subset $\Omega\subset\R^{d+1}$ bounded by an embedded, compact 
$C^1$-hypersurface $\Sigma=\partial\Omega$ is a Caccioppoli
set, satisfying $\Sigma=\Sigma^\ast $ and $\textnormal{Per}(\Omega)=\mathcal{H}^d(\Sigma)$. The measure-theoretic outward unit-normal coincides in this case with the classical continuous outward unit-normal defined on all of $\Sigma$; see \cite{FinitePerimeterSets} for more facts on Caccioppoli sets. 
The outward unit-normal $n_\Omega$ of a Caccioppoli set $\Omega$ can be used to express tangent-point energies in the co-dimension one case; similarly, e.g., as the unit-normal of a smooth hypersurface is used to define its fractional  Willmore energy  \cite{blatt2025fractional}.

\begin{defn}[Generalized tangent-point energies for hypersurfaces]
Let $d \in \mathbb{N}, p,q>0$ and $\Sigma=\partial \Omega$ be the boundary of a bounded Caccioppoli set $\Omega \subset \mathbb{R}^{d+1}$ with measure-theoretic outward unit-normal $n_\Omega: \Sigma^\ast\to\R^{d+1}$. Then the \emph{generalized tangent-point energies} are defined as
   \[
   \TP^{(p,q)}(\Sigma):= \int_{\Sigma ^\ast} \int_{\Sigma ^\ast} \frac{| \langle x-y,n_\Omega(x) \rangle |^q}{|x-y|^p} \, d \mathcal{H}^d(y) \, d \mathcal{H}^d(x) .
   \]
	If $p=2q$, this reduces to the classical \emph{tangent-point
		energies}
   \[
   \TP_q(\Sigma):=2^q \TP^{(2q,q)}(\Sigma)= 2^q
   \int_{\Sigma ^\ast} \int_{\Sigma ^\ast} \left(\frac{| \langle x-y,n_\Omega(x) \rangle |}{|x-y|^2} \right)^q \, d \mathcal{H}^d(y) \, d \mathcal{H}^d(x) .
   \]
\end{defn}

These energies are
non-local curvature energies  with self-avoidance properties. 
Their name originates from the fact that the Lagrangian of the classical tangent-point energies  $\TP_q$ is chosen to be the $q$-th inverse power of the \textit{tangent-point radius} $r_{\TP}(x,y)$, defined as the radius of the unique $d$-dimensional sphere through $x$ and $y$ that is tangent to $\Sigma$ at $x$. More concretely, we have the formula
\begin{equation}
	\label{eq: def tp-radius 1} 
	\frac{| \langle x-y,n_\Omega(x) \rangle |}{|x-y|^2}= \frac{1}{2 r_{\TP}(x,y)} \, .
\end{equation}
Classical tangent-point energies for sets of all dimensions and co-dimensions
were investigated by Strzelecki and the second author, who proved that
finite energy implies that the set is an embedded submanifold with an optimal
$C^{1,\alpha}$-regularity \cite[Theorem~1.3]{strzelecki2012tangent},
\cite[Theorem~1.4]{strzeleckivdm2}. Blatt \cite[Corollary~1.2]{blatt_2013b}
characterized the respective energy spaces of $\TP_q$ as particular Banach spaces, the fractional Sobolev-Slobodetski\u{\i} spaces $W^{2-(d/q),q}$ for $q>2d$. Motivated by regularity questions for critical knots, Blatt and Reiter admitted different powers $p$ and $q$ for numerator and denominator in the energy density for curves $\gamma$, to arrive at the generalized tangent-point energy $\TP^{(p,q)}(\gamma)$. The corresponding energy spaces $W^{(p-1)/q,q}$ 
\cite[Theorem~1.1]{blatt2012regularity} also include the one-parameter family of Hilbert spaces $W^{(p-1)/2,2}$. For $p\in (4,5)$ in the Hilbert case  $q=2$, they could show $C^\infty$-smoothness of $\TP^{(p,2)}$-critical curves \cite[Theorem~1.5]{blatt2012regularity}. The energy spaces of generalized
tangent-point energies $\TP^{(p,q)}$ on sets of arbitrary dimensions and co-dimensions were characterized by K\"afer in \cite[Theorem~5.40]{Kfer:820789} 
as the Sobolev-Slobodetski\u{\i} spaces $W^{(p-d)/q,q}$ if $q>d$ and $p\in [q+2d,2q+d)$.  In accordance with the case $d=1$ in \cite{blatt2012regularity}, we call the line $p=q+2d$ in the parameter range \emph{critical}, whereas the open region $p\in (q+2d,2q+d)$ is the \emph{subcritical} range.

\medskip 

{\bf Main Results.}
The goal of this article is to establish, for a broad range of parameter values $p$ and $q$ also beyond
the subcritical range, sharp lower bounds on the scale-invariant quotients 
\begin{equation}
\label{eq: scale invariant quotients}
    \TP^{(p,q)}(\Sigma) / \mathcal{H}^d(\Sigma ^\ast)^{\frac{2d+q-p}{d}} \quad \textnormal{and} \quad \TP^{(p,q)}(\Sigma) / \mathcal{L}^{d+1}(\Omega)^{\frac{2d+q-p}{d+1}} ,
\end{equation}
uniquely attained at round hyperspheres. 

To formulate our first main result for the classical tangent-point energy $\TP_q$ we denote by
\begin{equation}
	\label{eq: def r_d}
	r_{d}:= \left( \frac{\Gamma ( (d+1)/2)}{2 \pi ^{(d+1) /2}} \right)^{1/d}
\end{equation}
the radius of a $d$-dimensional round sphere $S\subset\R^{d+1}$ with $\mathcal{H}^{d}(S)=1$.
\begin{thm}[Classical tangent-point energies]
    \label{thm: main 1}
     Let $d \in \mathbb{N}$, $q \in [d, \infty)$ and $\Sigma=\partial \Omega$ be the boundary of a bounded Caccioppoli set $\Omega \subset \mathbb{R}^{d+1}$ with perimeter $\rm{Per}(\Omega)=\mathcal{H}^{d}(\Sigma^\ast)>0$. Then 
    \begin{equation}
    \label{eq: main theorem}
    \frac{\TP _{q} (\Sigma)}{\mathcal{H}^d(\Sigma ^\ast)^{\frac{2d-q}{d}}} \geq 
    \frac{\TP _{q} (S)}{\mathcal{H}^d(S)^{\frac{2d-q}{d}}}=\frac{1}{r_d^q} \, ,
    \end{equation}
    where $S$ denotes a $d$-dimensional sphere. If $q>d$, equality holds if and only if $\Sigma^\ast$ is a round sphere.
\end{thm}
Notice that for $q\leq 2d$, \Cref{thm: main 1} immediately implies for $\TP_q$ the analogous bound on the second quotient in \eqref{eq: scale invariant quotients}; see \Cref{rmk: implications of quotient bounds}.

Our second main theorem deals with generalized tangent-point
energies $\TP^{(p,q)}$ and covers a broader range of parameter values, which includes the classical case for $d \geq 3$. It provides a sharp lower bound, uniquely attained at round spheres, on both quotients in \eqref{eq: scale invariant quotients}. Remarkable is here the stronger lower bound on the second quotient \eqref{eq: TP^(p,q) minimum sets 2} in the subcritical case, which in particular constitutes an improvement of \Cref{thm: main 1} for $d \geq 3$ and $p=2q$.

\begin{thm}[Generalized tangent-point energies]
\label{thm: main 2}
    Let $d \in \mathbb{N}$, $q \in [3, \infty), p \in [\max \{q+d, 2q-4\},2q-3+d]$ and $\Sigma=\partial \Omega$ be the boundary of a bounded Caccioppoli set $\Omega \subset \mathbb{R}^{d+1}$ with perimeter $\rm{Per}(\Omega)=\mathcal{H}^{d}(\Sigma^\ast )>0$. Then
    \begin{equation}
    \label{eq: TP^(p,q) minimum sets 1}
    \frac{\TP^{(p,q)}(\Sigma)}{\mathcal{H}^d(\Sigma ^\ast)^{\frac{2d+q-p}{d}}}\geq \frac{\TP^{(p,q)}(S)}{\mathcal{H}^d(S)^{\frac{2d+q-p}{d}}}
    \end{equation}
    and 
    \begin{equation}
     \label{eq: TP^(p,q) minimum sets 2}
    \frac{\TP^{(p,q)}(\Sigma)}{\mathcal{L}^{d+1}(\Omega)^{\frac{2d+q-p}{d+1}}}\geq \frac{\TP^{(p,q)}(S)}{\mathcal{L}^{d+1}(B)^{\frac{2d+q-p}{d+1}}} \, ,
    \end{equation}
    where $S=\partial B$ denotes a $d$-dimensional round sphere and $B$ the ball enclosed by $S$. If $p>q+d$, equality in either inequality holds if and only if $\Sigma^\ast$ is a round sphere. 
\end{thm}

As mentioned above, an important subclass of Caccioppoli sets is that of open bounded domains $\Omega \subset \mathbb{R}^{d+1}$ with $C^1$-boundary $\Sigma=\partial\Omega =\partial ^\ast \Omega$. Since these sets are of particular interest in geometric knot theory and low-dimensional topology, we highlight the following immediate consequence of the two main theorems.
\begin{cor}
\label{cor: TP^(p,q) min C^1 boundaries}
    Let $d \in \mathbb{N}$, $p, q>0$ and $\Sigma=\partial \Omega$ be the 
    $C^1$-boundary of a non-empty open bounded set $\Omega \subset \mathbb{R}^{d+1}$.
    \begin{enumerate}[\rm (i)]
        \item If $q \in [d , \infty)$, then \eqref{eq: main theorem} holds true. If $q >d$, then equality holds in \eqref{eq: main theorem} if and only if $\Sigma$ is a round sphere.
        \item If $q \in [3, \infty), p \in [\max \{q+d, 2q-4\},2q-3+d]$, then \eqref{eq: TP^(p,q) minimum sets 1} and \eqref{eq: TP^(p,q) minimum sets 2} hold true. If $p>q+d$, then equality holds in \eqref{eq: TP^(p,q) minimum sets 1} or \eqref{eq: TP^(p,q) minimum sets 2} if and only if $\Sigma$ is a round sphere.
    \end{enumerate}
\end{cor} 

Before indicating the crucial steps of the proofs, let us recall the geometric meaning of the tangent-point radius $r_\textnormal{TP}$ in \eqref{eq: def tp-radius 1} defining the energy density of $\TP_q$. This helps us to
characterize round spheres as the only boundaries of convex sets for which the integrand of $\TP_q$ is constant; see \Cref{lem: sphere iff constant tangent point radius rectifiable} and \Cref{cor: sphere constant integrand}. For several geometric functionals on closed hypersurfaces, such rigidity of the energy density often comes together with sharp lower bounds, uniquely attained at spheres. Examples of 
such functionals uniquely minimized by spheres include the $L^p$-norms of the Gauss-Kronecker curvature for $p>1$, which  can be seen by combining the Fenchel-type theorem of \cite[Theorem~1]{GaussBoundGeneralDim} for the total Gauss-Kronecker curvature with Jensen's inequality and the rigidity of the pointwise Gauss-Kronecker curvature proved in \cite[Theorem~2]{GaussUniqueGeneralDim}. Moreover, the $p$-Willmore energies for 
$p\geq d\geq 2$ exhibit such spherical rigidity as well \cite[Theorem~3]{WillmoreGeneralDim}.

If one replaces the inner integral of the classical tangent-point energy $\TP_q$ by a supremum, one obtains as a non-differentiable functional the $L^q$-norm of the
\emph{global tangent-point curvature}
\begin{equation}\label{eq:global-tp}
	\mathcal{K}[\Sigma](x):=\sup_{y\in\Sigma}\frac1{r_\TP(x,y)}\quad\textnormal{for $x\in
		\Sigma,$}
\end{equation}
which provides much more geometric rigidity because the supremum serves as a hard-core repulsive potential for the global interaction between different surface points. Indeed, using elementary geometry in combination with the classical isoperimetric inequality it was shown in \cite[Theorem~1.5]{vdMStrKol} that round hyperspheres uniquely minimize the $L^q$-norm of the global tangent-point curvature, that is, $\|\mathcal{K}[\Sigma]\|_{L^q(\Sigma,\mathcal{H}^d)}$, among a broad class of admissible sets $\Sigma\subset\R^{d+1}$ including closed $C^1$-hypersurfaces\footnote{In the context of the present paper, the rigidity of $\|\mathcal{K}[\Sigma]\|_{L^q(\Sigma,\mathcal{H}^d)}$ is an immediate consequence of \Cref{thm: main 1}}. This resembles an earlier rigidity result on the $L^q$-norm of a global curvature function on curves; see \cite[Lemma 7]{strzelecki-vdm_2007}. 

Thus, one may expect the same to hold true for the classical tangent-point energies $\TP_q$. But we should emphasize that, in contrast to these single-integral curvature energies, tangent-point energies are \emph{non-local} double integrals with interactions between  points that are intrinsically far apart on $\Sigma$. This makes the analysis harder.

On the other hand, $\TP^{(p,q)}$ is the higher-dimensional analogue of the generalized tangent-point energies on curves, for which D\"ohrer and the first author show in \cite[Theorems~1.1~and~1.2]{döhrer2026fencheltheoremgaussmaps} that
the round circle uniquely minimizes $\TP^{(p,q)}$ among all closed curves under a fixed-length constraint for a wide range of parameters $p$ and $q$, extending an earlier rigidity result of Blatt and Volkmann on the classical tangent-point energy $\TP_q$ for curves in \cite[Proposition~5.10~\&~Corollary~5.12]{volkmann1}. In light of that, \hyperref[thm: main 1]{Theorems \ref*{thm: main 1}} and
\hyperref[thm: main 2]{\ref*{thm: main 2}} are natural generalizations to
arbitrary dimensions, including a substantial range of decoupled parameters
$p$ and $q$, where the integrand is \emph{not} constant on spheres.
We should point out, however, that not all tools developed in \cite{döhrer2026fencheltheoremgaussmaps} are 
applicable in the higher-dimensional case. In particular, the first part 
of the central \textit{Fenchel Theorem for Gauss Maps} \cite[Theorem~1.4]{döhrer2026fencheltheoremgaussmaps} seems to lack a 
higher-dimensional analogue. Even the earlier inequality on the classical tangent-point energies on curves does not appear to be generalizable to higher dimensions by a mere adaptation of the arguments in \cite{volkmann1}. 

\medskip

{\bf Strategy.}
The proofs of \hyperref[thm: main 1]{Theorems \ref*{thm: main 1}} and
\hyperref[thm: main 2]{\ref*{thm: main 2}} combine elementary inequalities with integral geometric formulas, such as a generalization of Bushling's identity \cite[Theorem, p.~3]{bushling2024singular} 
\begin{equation}
\label{eq: Bushling formula}
    \frac{1}{\alpha_{d}}\int_{\Upsilon} \int_{\Upsilon} \frac{\langle x-y, n(y) \rangle \langle y-x, n(x) \rangle}{|x-y|^{d+2}} 
    \, d \mathcal{H}^{d}(y) \, d \mathcal{H}^{d}(x)= \mathcal{H}^{d}(\Upsilon)
\end{equation}
for any $d$-rectifiable set $\Upsilon \subset \mathbb{R}^{d+1}$ that possesses a unit-normal vector field $n$ satisfying the so-called Orientation Cancellation Condition (OCC). Here, $\alpha_d$ is given by the constant $\frac{1}{2} \int_{\mathbb{S}^{d}} |\langle \omega, e_1 \rangle | \, d \mathcal{H}^{d}(\omega)$. Roughly speaking, the OCC means that the sum of all signed intersections between $\Upsilon$ and a given oriented line $\ell$ through a base point $x \in \Upsilon$ is zero (up to $\mathcal{H}^{d}$-null sets). Bushling showed that the OCC is satisfied by the reduced boundary $\Sigma ^\ast=\partial ^\ast \Omega$ of any bounded Caccioppoli set and its measure theoretic outward unit-normal $n_\Omega$ \cite[Lemma~2]{bushling2024singular}. In fact, Bushling's proof implies the stronger statement that for such a line $\ell$ and a point $x \in \Sigma ^\ast$ there is an even number of intersection points with alternating signs (see \Cref{lem: finite perimeter sets satisfy SOCC} of the present paper). For the proof of \Cref{thm: main 2} we will invoke the fact that $(\Sigma ^\ast, n_\Omega)$ satisfies this stronger property.

A central object in the proof of integral identities such as \eqref{eq: Bushling formula} is what we refer to as the \textit{Gauss map} of the reduced boundary $\Sigma^\ast$ for fixed $x \in \Sigma ^\ast$
\[
\phi_x: \Sigma ^\ast \setminus \{x\} \rightarrow \mathbb{S}^d, \quad  \phi_x(y)=\frac{y-x}{|y-x|}  .
\]
Via the coarea formula this map links the left-hand side of \eqref{eq: Bushling formula} with averages over signed intersections between $\Sigma ^\ast$ and lines in $\mathbb{R}^{d+1}$, which can then be analyzed using the OCC.

From \eqref{eq: Bushling formula} one immediately obtains the inequality
\begin{equation}
\label{eq: Steinerberger/Bushling inequality}
    \frac{1}{\alpha_{d}}\int_{\Sigma^\ast} \int_{\Sigma^\ast} \frac{|\langle x-y, n_\Omega(y) \rangle \langle y-x, n_\Omega(x) \rangle|}{|x-y|^{d+2}} 
    \, d \mathcal{H}^{d}(y) \, d \mathcal{H}^{d}(x) \geq  \mathcal{H}^{d}(\Sigma ^\ast)
\end{equation}
for the reduced boundary $\Sigma ^\ast$ of any Caccioppoli set \cite[Corollary, p.~4]{bushling2024singular}, which was previously shown to hold for $C^1$-boundaries by Steinerberger in \cite{steinerberger2022inequality}. In both works, equality was shown to hold if and only if $\Omega$ is convex (up to $\mathcal{L}^{d+1}$-equivalence). For our approach, we will extend in \Cref{thm: bushling theorem general} the perimeter-identity \eqref{eq: Bushling formula} by generalizing the integrand of the left-hand side, to obtain a generalization of inequality \eqref{eq: Steinerberger/Bushling inequality} stated in \Cref{cor: lower bound for F_G for sets}. As a special case we will obtain for $\varrho \geq 1$ the inequality
\begin{equation}
	\label{eq: F_q inequality}
  \frac{1}{b_{d, \varrho}} \int_{\Sigma ^\ast }\int_{\Sigma ^\ast}
\left(\frac{| \langle x-y,n_\Omega(x) \rangle |}{|x-y|} \right)^{\varrho-1} \,
\frac{| \langle x-y,n_\Omega(y) \rangle |}{|x-y|^{d+1}} 
\, d \mathcal{H}^{d}(y) \, d \mathcal{H}^{d}(x) \geq \mathcal{H}^{d}(\Sigma ^\ast) \, ,  
\end{equation}
where $b_{d, \varrho}= \frac{1}{2} \int_{\mathbb{S}^{d}} |\langle \omega, e_1 \rangle |^{\varrho -1} \, d \mathcal{H}^{d}(\omega)$ and equality holds here exactly for boundaries of convex sets; see \Cref{cor: lower bound on F_q for sets}. Inequality \eqref{eq: F_q inequality} is not only a straightforward   extension of \eqref{eq: Steinerberger/Bushling inequality} but may also be
interpreted as a higher-dimensional analogue of  \cite[Lemma~5.12]{döhrer2026fencheltheoremgaussmaps} albeit in co-dimension one.

It turns out that the left-hand side of \eqref{eq: F_q inequality} with $\varrho=q$ is a sharp lower bound of $\frac{1}{b_{d, q}} \TP^{(q+d,q)}$, which implies
\begin{equation}
\label{eq: TP^(q+d,q) inequality}
	 \frac{\TP^{(q+d,q)}(\Sigma)}{\mathcal{H}^{d}(\Sigma ^\ast)} \geq 
	\frac{\TP^{(q+d,q)}(S)}{\mathcal{H}^d(S)} .
\end{equation} 
We refer to the special line $p=q+d$ in the parameter range as the \textit{lower limit case}, in accordance with the terminology used in \cite{döhrer2026fencheltheoremgaussmaps} for $d=1$. This terminology is justified by the fact that for every $p<q+d$ the infimum of the first quotient in \eqref{eq: scale invariant quotients} among all boundaries of bounded Caccioppoli sets is equal to zero. Indeed, a minimizing sequence is, e.g., given by the boundaries of the increasingly flat boxes $\Omega_k:=(0,1)^d \times (0, \frac{1}{k})$ as $k \to \infty$. Furthermore, notice that the lower limit case has no overlap with the critical case $p=q+2d$ or subcritical case $p>q+2d$. In fact, the integrand of $\TP^{(q+d,q)}$ is not singular enough to provide sufficiently strong self-avoidance properties (see \Cref{cor: upper bound on TP}).
Due to the constancy of the integrand of $\TP_q$ on spheres as shown in \Cref{cor: sphere constant integrand}, inequality \eqref{eq: main theorem} follows from \eqref{eq: TP^(q+d,q) inequality} with the choice $q=d$. 

For the generalized tangent-point energies $\TP^{(p,q)}$, however, we lose this rigidity of the Lagrangian for decoupled parameters $p$ and $q$. So, we need to control an additional term, the \textit{signed symmetrized tangent-point functional} $\TP_{\textnormal{ssym}}^{(\alpha,2)}$ of $\Sigma =\partial \Omega$, defined as the iterated integral
\[
\TP_{\textnormal{ssym}}^{(\alpha,2)}(\Sigma)
:= 
\int_{\Sigma ^\ast} \int_{ \Sigma ^\ast}
\frac{ \langle x-y , n_\Omega(x) \rangle \, \langle y-x,n_\Omega(y) \rangle }{|x-y|^\alpha} 
d \mathcal{H}^{d}(y) \, d\mathcal{H}^{d}(x) \, , \quad \alpha \in [0, d+2].
\]
For $\alpha \in (0, d+1)$ the term $\TP_{\textnormal{ssym}}^{(\alpha,2)}(\Sigma)$ can be related to the \textit{Riesz energy} of $\Omega$, given by
\[
R_\alpha(\Omega):= \int_{\Omega} \int_{\Omega} |x-y|^{-\alpha} 
\, dx \, dy \, , \quad \alpha \in (0,d+1) .
\]
More precisely, with the Gauss-Green formula we will establish in \Cref{lem: TP_sym and Riesz relation} the identity  
\begin{equation}
\label{eq: TP_ssym Riesz identity}
\TP_{\textnormal{ssym}}^{(\alpha,2)}(\Sigma)= (-\alpha+d+2)(-\alpha+d+1) \, R_\alpha (\Omega) \, , \quad \alpha \in (0, d+1) \, ,
\end{equation}
which extends to 
\[
\TP_{\textnormal{ssym}}^{(0,2)}(\Sigma)=(d+2) (d+1) \, \mathcal{L}^{d+1}(\Omega)^2 
\quad \textnormal{and} \quad 
\TP_{\textnormal{ssym}}^{(d+1,2)}(\Sigma)= \mathcal{H}^{d}(\mathbb{S}^{d}) \, \mathcal{L}^{d+1}(\Omega) .
\]
The Riesz energy is uniquely maximized by balls among all Lebesgue-measurable sets $\Omega \subset \R^{d+1}$ of fixed volume \cite[Equation~(1.2)]{RieszIso}, which is a consequence of the classical Riesz rearrangement inequality.

While several equivalent representations, including double boundary integrals, are known for $R_\alpha$, the identity \eqref{eq: TP_ssym Riesz identity} is to the best of the authors' knowledge not yet established in the literature. However, it is inspired by several identities established in a series of papers by O'Hara and Solanes \cite{o2012centers,o2015planar,o2018regularized}. Among these, a similar double boundary integral representation of $R_\alpha$ is given by \cite[Lemma~4.1]{o2018regularized} alongside an analogous one for the more singular \textit{fractional perimeter} $P_\alpha$, $\alpha \in (d+1,d+2)$. Similarly, \eqref{eq: TP_ssym Riesz identity} can be extended to an identity between $\TP_{\textnormal{ssym}}^{(\alpha,2)}$ and $P_\alpha$ for $\alpha \in (d+1,d+2)$. Furthermore, for $d=1$, \eqref{eq: TP_ssym Riesz identity} can be viewed as a less singular version of the identity in \cite[Proposition~3.13]{o2015planar}. 

By means of \eqref{eq: TP_ssym Riesz identity} we will establish a sharp lower bound on $R_\alpha$ (see \Cref{lem: sharp lower bound on Riesz energy}) in terms of the chord power integral 
\begin{equation}
	\label{eq: def I_alpha}
I_\alpha(\Omega):= \int_{\Sigma ^\ast} \int_{S_x} 
l(x, \omega) ^{-\alpha+d+2} \, \langle -\omega , n_\Omega(x) \rangle
\, d \mathcal{H}^{d}(\omega) \, d \mathcal{H}^{d}(x) \, , \quad \alpha \in [0,d+1],
\end{equation}
where $S_x:= \{ \omega \in \mathbb{S}^{d} \, | \, \langle \omega , n_\Omega (x) \rangle <0 \}$ denotes the set of inward directions at $x \in \Sigma ^\ast$, and 
\[
l(x,\omega):= \sup \{t \geq0  : \, \chi _{\Omega} (x+s \omega) =1 \quad \textnormal{for a.e. } s \in (0,t) \} 
\]
denotes the length of the chord between $x$ and the point at which the ray $x + \R ^+ \omega$ first exits $\Omega$. More precisely, we will show that for any bounded Caccioppoli set $\Omega$ there holds
\begin{equation}
	\label{eq: sharp lower bound on R_alpha via I_alpha}
	(-\alpha+d+2)(-\alpha+d+1) \, R_\alpha (\Omega) \geq I_\alpha(\Omega) , \quad \alpha \in (0,d+1)
\end{equation}
and
\[
(d+2) (d+1) \, \mathcal{L}^{d+1}(\Omega)^2 \geq I_{0}(\Omega)
\]
with equality if and only if $\partial ^\ast \Omega$ is the reduced boundary of a convex set, as well as
\[
\mathcal{H}^{d}(\mathbb{S}^{d}) \, \mathcal{L}^{d+1}(\Omega) = I_{d+1}(\Omega).
\]
Chord power integrals such as $I_\alpha$, their more singular extensions, and their relations to the Riesz energy and the fractional perimeter have been studied in \cite{ChordPowerIntegral,LittlePrince,FracPerimeterUpperBound}. 

\medskip

{\bf Structure of the paper.}
In \Cref{sec: Preliminaries} we recall some basic definitions and results in the realm of Caccioppoli sets mostly based on \cite{FinitePerimeterSets}, but we note that \cite{BoundedVariation} also contains a lot of nuanced results on that matter. The machinery from \cite{bushling2024singular} and its consequences for our situation is presented as well, before we conclude \Cref{sec: Preliminaries} with two statements concerning sets with constant tangent-point radius. These will be used later on to ensure that certain lower bounds for $\TP^{(p,q)}$ are attained uniquely at spheres.

\Cref{sec: Proof of the main theorem} is devoted to the proof of \Cref{thm: main 1}, consisting of two steps. In the first step, we prove a generalization of \eqref{eq: Bushling formula}, from which we conclude \eqref{eq: F_q inequality} with equality exactly for boundaries of convex sets. From that, we deduce in the second step first \eqref{eq: TP^(q+d,q) inequality}, and subsequently \eqref{eq: main theorem}. With the help of \Cref{lem: sphere iff constant tangent point radius rectifiable} we prove rigidity of spheres in \eqref{eq: main theorem} for $q>d$, concluding the proof of \Cref{thm: main 1}. 

In \Cref{sec: More parameters} we prove \Cref{thm: main 2}. To this end, we first establish in \Cref{lem: TP_sym and Riesz relation} the identity \eqref{eq: TP_ssym Riesz identity}, from which we deduce in \Cref{lem: sharp lower bound on Riesz energy} inequality \eqref{eq: sharp lower bound on R_alpha via I_alpha} and its extensions to $\alpha=0$ and $\alpha=d+1$. From that we derive a sharp lower bound on both quotients in \eqref{eq: scale invariant quotients}, involving the perimeter and the term $-R_{p-(2q-4)}$.
\Cref{prop: Riesz maximum} demonstrates how to translate the sharp lower bound to an expression that only contains the perimeter and the volume of $\Omega$. Finally, we conclude the proof of \Cref{thm: main 2} including its rigidity discussion via the classical isoperimetric inequality.

\section{Preliminaries}
\label{sec: Preliminaries}

The first definition and result is formulated for $k$-dimensional subsets of $\R^m$. As soon as we focus on the co-dimension one case we specify to $m:=d+1$.

\begin{defn}[Rectifiable sets]
\label{defn: rectifiable sets}
A set $\Sigma \subseteq \mathbb{R}^{m}$ is \textit{k-rectifiable} if $\mathcal{H}^{k}(\Sigma) < \infty$ and there exists a countable family $\{F_i\}_{i \in I}$ of Lipschitz maps $F_i: A_i \subseteq \mathbb{R}^k \rightarrow \mathbb{R}^m$ and an $\mathcal{H}^k$-null set $\Sigma_0 \subseteq \Sigma$ such that 
\[
\Sigma = \Sigma_0 \cup \bigcup _{i \in I} F_i (A_i).
\]
If $\Sigma$ is $k$-rectifiable, then, at $\mathcal{H}^k$-a.e. $y \in \Sigma$, there exists a unique \textit{approximate tangent space} $T_y \Sigma \subset \mathbb{R}^m$, which coincides with the classical tangent space when $\Sigma$ is a smooth hypersurface. We call $n: \Sigma \rightarrow \mathbb{S}^{m-1}$ a \textit{measurable unit-normal vector field} on $\Sigma$ if it is $\mathcal{H}^k$-measurable and $n(y)$ is orthogonal to $T_y\Sigma$ for $\mathcal{H}^{k}$-a.e. $y \in \Sigma$ at which the approximate tangent space is uniquely defined.
\end{defn}

We will need the following coarea formula for Lipschitz maps between equal-dimensional rectifiable sets, which can be understood as a special case of \cite[Theorem~3.2.22]{Federer}. To this end, we denote by $|Jf|$ the tangential Jacobian of a (locally) Lipschitz map $f:U \rightarrow V$ between $k$-rectifiable sets $U,V \subset \mathbb{R}^{m}$, which coincides with the $k$-dimensional Jacobian $J_kf$ in \cite[Definition~3.2.1]{Federer}. The cited theorem is formulated for globally Lipschitz maps $f$, yet under suitable integrability or nonnegativity assumptions extends immediately to locally Lipschitz $f$ via a partition of $U$ into rectifiable sets on which $f$ is globally Lipschitz.

\begin{prop}[Coarea Formula for equal-dimensional sets]
	\label{prop: coarea formula}
	Let $U,V \subset \mathbb{R}^{m}$ be $k$-rectifiable sets, and let $f: U \rightarrow V$ be locally Lipschitz. Then the tangential Jacobian $|Jf(x)|$ exists for $\mathcal{H}^k$-a.e. $x \in U$ and the following identity holds for every nonnegative measurable function $g$, with both sides interpreted as extended nonnegative integrals:
	\[
	\int_U g(x) \, |J f(x) | \, d \mathcal{H}^k (x) = \int_{V} \sum _{z \in f^{-1}(y)} g(z) \, d \mathcal{H}^{k}(y)  .
	\]
	The same identity holds for signed measurable functions $g$ satisfying 
	\[
	\int_U |g(x)| \, |J f(x) | \, d \mathcal{H}^k (x) < \infty  .
	\]
\end{prop}

\begin{prop}[De Giorgi's structure theorem]
	\label{prop: Gauss-Green}
	If $\Omega \subseteq \mathbb{R}^{d+1}$ is a Caccioppoli set, then the generalized Gauss-Green formula holds true:
	\[
	\int_{\Omega} \textnormal{div}(\varphi) = \int_{\partial ^\ast \Omega} \langle \varphi, n_\Omega \rangle \, d \mathcal{H}^{d}  \quad \forall \varphi \in C^1_c(\mathbb{R}^{d+1}, \mathbb{R}^{d+1})  .
	\]
\end{prop}
\begin{proof}
	The formula follows directly from \cite[Theorem~15.9]{FinitePerimeterSets} applied to each component of $\varphi$ after taking the trace on both sides of the resulting matrix equation.
\end{proof}

We call a hyperplane $P \subset \mathbb{R}^{d+1}$ a supporting hyperplane of a set $\Omega \subset \mathbb{R}^{d+1}$ at $x \in \partial \Omega$ if $x \in P$ and $\Omega$ is entirely contained in one of the two closed half-spaces bounded by $P$. The following proposition constitutes a \textit{Supporting Hyperplane Theorem} for Caccioppoli sets with an explicit formula involving the measure-theoretic outward unit-normal $n_\Omega$.

\begin{prop}[Supporting hyperplane theorem]
	\label{prop: Supporting Hyperplane Theorem strong}
	Let $\Omega \subset \mathbb{R}^{d+1}$ be a set with its open interior denoted by $\Omega^{\mathrm{o}}$. 
	\begin{enumerate}[\rm (i)]
		\item If $\Omega$ is bounded and convex, then it is a Caccioppoli set and for every $x \in \partial \Omega$ there exists a supporting hyperplane. Moreover, for every $x \in \partial ^\ast \Omega$ a supporting hyperplane is given by $x+n_\Omega (x) ^\perp$ and 
		\[
		\Omega \subset \{ y \in \mathbb{R}^{d+1} \, | \, \langle y-x, n_\Omega (x) \rangle \leq 0 \}  .
		\]
		If in addition $\mathcal{H}^{d}(\partial ^\ast \Omega)>0$, then $\mathcal{L}^{d+1}(\Omega)>0$ and $\mathcal{H}^{d}(\partial ^\ast \Omega)=\mathcal{H}^{d}(\partial \Omega)$.
		\item Suppose $\Omega$ is a closed Caccioppoli set with non-empty interior $\Omega^{\mathrm{o}}$. If $\Omega$ has a supporting hyperplane at every $x \in \partial ^\ast \Omega$, then $\Omega$ is convex.
	\end{enumerate}
\end{prop}

\begin{proof}
	According to \cite[Example~12.6]{FinitePerimeterSets} bounded convex sets are Caccioppoli sets. The claim in (i) about the possible candidates for supporting hyperplanes together with the succeeding formula can be found at \cite[p.~2038]{FinitePerimeterSetsConvex2}. Suppose now that $\Omega$ is bounded and convex. If $\mathcal{H}^{d}(\partial ^\ast \Omega)>0$, then by \cite[Remark~12.4]{FinitePerimeterSets} $\mathcal{L}^{d+1}(\Omega)>0$ and due to convexity of $\Omega$ also $\mathcal{L}^{d+1}(\Omega^{\mathrm{o}})>0$. Thus, $\overline{\Omega^{\mathrm{o}}}=\overline{\Omega}$ or equivalently $\partial \Omega^{\mathrm{o}}= \partial \Omega$. Since $\mathcal{L}^{d+1}(\partial \Omega)=0$, we also have $\partial ^\ast \Omega=\partial ^\ast \Omega^{\mathrm{o}}$. Therefore we conclude by \cite[Example~12.6]{FinitePerimeterSets}
	\[
	\mathcal{H}^{d}(\partial^\ast \Omega)= \mathcal{H}^{d}(\partial ^\ast \Omega^{\mathrm{o}})=\rm{Per}(\Omega^{\mathrm{o}})=\mathcal{H}^{d}(\partial \Omega^{\mathrm{o}})=\mathcal{H}^{d}(\partial \Omega)  .
	\]
	Part $(ii)$ follows directly from \cite[Theorem~1]{FinitePerimeterSetsConvex}.
\end{proof}

\begin{nota}
	\label{nota: L and l}
	For each $x \in \mathbb{R}^{d+1}$ and $\omega \in \mathbb{S}^{d}$, let	
		\[	
		h_{x,\omega}: \mathbb{R} \rightarrow \mathbb{R}^{d+1}, \quad t \mapsto x+ t \omega
		\quad \textnormal{and} \quad 
		L_{x,\omega}:=h_{x,\omega}(\mathbb{R}) =x+ \textnormal{span } \omega  .
		\]
	Define further for a measurable set $\Omega \subset \mathbb{R}^{d+1}$
	\[
	l(x,\omega):= \sup \{t \geq0  : \, \chi _{\Omega} (x+s \omega) =1 \quad \textnormal{for a.e. } s \in (0,t) \}  .
	\]
\end{nota}

\begin{defn}[Orientation cancellation condition (OCC) (\cite{bushling2024singular}, Definition, p.~5)]
\label{defn: OCC} 
Given a $d$-rectifiable set $\Sigma \subset \mathbb{R}^{d+1}$ and a measurable unit-normal vector field $n: \Sigma \rightarrow \mathbb{S}^{d}$, we say the \textit{orientation cancellation condition} (or \textit{OCC}) \textit{is satisfied} at a point $x \in \Sigma$ if the equation
\[
\sum_{y \in \Sigma \cap L_{x,\omega}} \textnormal{sgn } \langle \omega, n(y) \rangle =0
\]
holds for $\mathcal{H}^{d}$-a.e. $\omega \in \mathbb{S}^{d}$, where sgn is the signum function with sgn $0:=0$.
If the OCC is satisfied at $\mathcal{H}^{d}$-a.e. $x \in \Sigma$, we say that $(\Sigma, n)$ \textit{satisfies the orientation cancellation condition}.
\end{defn}

By \cite[Lemma~2]{bushling2024singular} reduced boundaries of Caccioppoli sets satisfy the OCC as a consequence of the following stronger statement.
\begin{lem}
\label{lem: finite perimeter sets satisfy SOCC}
    Let $\Omega \subset \mathbb{R}^{d+1}$ be a bounded Caccioppoli set. Then for $\mathcal{H}^{d}$-a.e. $x \in \partial ^\ast \Omega$ and $\mathcal{H}^{d}$-a.e. $\omega \in \mathbb{S}^{d}$ there exist $k \in \mathbb{N}_0$, $t_1<t_2<...<t_{2k}$, such that $h_{x,\omega}^{-1}(\Omega)$ is $\mathcal{L}^1$-equivalent to $\bigcup_{i=1}^{k} (t_{2i-1}, t_{2i})$. Moreover, $L_{x, \omega} \cap \partial ^\ast \Omega=\{x+t_i \omega \, | \, 1\leq i\leq 2k\}$ and
    \[
    \textnormal{sgn } \langle \omega , n (x+ t_i \omega) \rangle = (-1)^i .
    \]
    In particular, $(\partial ^\ast \Omega, n_\Omega)$ satisfies the OCC.
\end{lem}

\begin{proof}
	We closely follow the proof of \cite[Lemma~2]{bushling2024singular} with the identification $E=\Omega$ and $v_E=n_\Omega$.
    According to \cite[Theorem, p.~9]{bushling2024singular} for $\mathcal{H}^{d}$-a.e. $x \in \partial ^\ast \Omega$ and $\mathcal{H}^{d}$-a.e. $\omega \in \mathbb{S}^{d}$ the set of discontinuities of the function $(\chi_\Omega^\ast)_x^\omega$ defined right before \cite[Equation~(3.2)]{bushling2024singular} is given by the set $(J_\Omega)_x^\omega:=\{t \in \R : x+ t \omega \in \partial ^\ast \Omega\}$, which by the succeeding discussion satisfies $\mathcal{H}^0((J_\Omega)_x^\omega) \in 2 \mathbb{Z}$. Furthermore $(\chi_\Omega^\ast)_x^\omega(t)=\chi_\Omega (x+t \omega)$ for $\mathcal{L}^1$-a.e. $t \in \R$. Fixing such $x$ and $\omega$, we obtain $t_1<t_2<...<t_{2k}$ by enumerating $(J_\Omega)_x^\omega$. Immediately, we have $L_{x, \omega} \cap \partial ^\ast \Omega=\{x+t_i \omega \, | \, 1\leq i\leq 2k\}$. The proof of \cite[Lemma~2]{bushling2024singular} continues to show that $(\chi_\Omega^\ast)_x^\omega(t)=0$ for all $t \in (t_i,t_{i+1})$, whenever $\langle \omega , n_\Omega (x+t_i \omega) \rangle >0$, leading to $\langle \omega, n_\Omega(x+t_{i+1} \omega) \rangle <0$. An argument analogous to the one in that proof shows that $(\chi_\Omega^\ast)_x^\omega(t)=1$ for all $t \in (t_i,t_{i+1})$, whenever $\langle \omega , n_\Omega (x+t_i \omega) \rangle <0$ and subsequently $\langle \omega, n_\Omega(x+t_{i+1} \omega) \rangle >0$, which proves the alternating sign claim of \Cref{lem: finite perimeter sets satisfy SOCC}. Finally, by $\mathcal{L}^{1}$-equivalence of $(\chi_\Omega^\ast)_x^\omega$ and $\chi_\Omega (x+ \cdot \, \omega)$ the previous discussion implies that $h_{x,\omega}^{-1}(\Omega)$ is $\mathcal{L}^1$-equivalent to $\bigcup_{i=1}^{k} (t_{2i-1}, t_{2i})$.
\end{proof}

Given a set $M \subset \mathbb{R}^{d+1}$ and a point $x \in M$, we consider the \textit{Gauss map} 
\begin{equation}
	\label{eq: def: Gauss map}
	\phi_x : M \setminus \{x \} \rightarrow \mathbb{S}^{d} , \quad \phi_x(y)= \frac{y-x}{|y-x|}  .
\end{equation}

As a consequence of \Cref{lem: finite perimeter sets satisfy SOCC} we obtain the following surjectivity property of the Gauss map.
\begin{lem}
\label{lem: surjectivity of Gauss map}
    Let $\Omega \subset \mathbb{R}^{d+1}$ be a bounded Caccioppoli set. For $z \in \partial ^\ast \Omega$ let $S_z:=\{\omega \in \mathbb{S}^{d} \, | \, \langle \omega , n_\Omega(z) \rangle <0 \}$. Then for $\mathcal{H}^{d}$-a.e. $x \in \partial ^\ast \Omega$, there holds that $S_x \subseteq  \phi_x (\partial ^\ast \Omega \setminus \{x\})$ up to a set of $\mathcal{H}^{d}$-measure zero. If $\Omega$ is convex, then $\mathcal{H}^{0}(\phi_x^{-1}(\omega))\leq 1$  for $\mathcal{H}^{d}$-a.e. $x \in \partial ^\ast \Omega$ and $\omega \in \mathbb{S}^{d}$.
\end{lem}

\begin{proof}
    By \Cref{lem: finite perimeter sets satisfy SOCC} we obtain for $\mathcal{H}^{d}$-a.e. $x \in \partial ^\ast \Omega$ and $\mathcal{H}^{d}$-a.e. $\omega \in S_x \subset \mathbb{S}^{d}$ a finite family $t_1<t_2<...<t_{2k}$ such that $L_{x, \omega} \cap \partial ^\ast \Omega=\{ x+ t_i \omega \, | \, 1 \leq i \leq 2k \}$ and $\langle \omega , n (x+ t_{2k} \omega) \rangle >0$.
    Fix such $x$ and $\omega $. Clearly, there exists $i^\ast \in \{1, ..., 2k\}$ such that $t_{i^\ast}=0$. Since $\langle \omega, n(x+t_{i^\ast})\rangle <0$, we conclude $t_{2k}>0$ and thus for $y=x+t_{2k} \omega \in \partial ^\ast \Omega \setminus \{x\}$ we have 
    \[
    \phi_x(y)= \frac{x+t_{2k} \omega - x}{|x+t_{2k} \omega -x|}= \textnormal{sgn}(t_{2k}) \omega=\omega .
    \]
    Thus, $\omega \in \phi_x (\partial ^\ast \Omega \setminus \{x\})$. Since this is true for $\mathcal{H}^{d}$-a.e. $x \in \partial ^\ast \Omega$ and $\mathcal{H}^{d}$-a.e. $\omega \in S_x$, the first claim follows.
    If $\Omega$ is convex, then by the proof of \cite[Corollary, p.~4]{bushling2024singular} for $\mathcal{H}^{d}$-a.e. $x \in \partial ^\ast \Omega$ and almost every $\omega \in \mathbb{S}^{d}$ there holds $\mathcal{H}^{0}(\partial ^\ast \Omega \cap L_{x, \omega})\leq 2$, which implies $\mathcal{H}^{0}(\phi_x^{-1}(\omega))\leq 1$.
\end{proof}

Combining \Cref{lem: finite perimeter sets satisfy SOCC} with the characterizing property of convex sets, that almost every line intersects its boundary in exactly two points, leads to the following lemma; see \cite[Inequality~(3.3)]{bushling2024singular}.

\begin{lem}
\label{lem: convex set nonlocal characterization}
    Let $\Omega \subset \mathbb{R}^{d+1}$ be a bounded Caccioppoli set. Then $\partial ^\ast \Omega$ is the reduced boundary of a bounded convex set if and only if
    \[
    \langle x-y , n_\Omega (y) \rangle \langle y-x, n_\Omega (x) \rangle \geq 0 \quad \textnormal{for } \mathcal{H}^{2d}-\textnormal{a.e. } (x,y) \in \partial ^\ast \Omega \times \partial ^\ast \Omega .
    \]
\end{lem}

In the following, let $S_{a,R}$ and $B_{a,R}$ denote, respectively, the $d$-dimensional sphere and the open $(d+1)$-dimensional ball in $\mathbb{R}^{d+1}$ with center $a \in \mathbb{R}^{d+1}$ and radius $R>0$.

\begin{lem}
\label{lem: formula for tp-radius}
    For every $x \in \mathbb{R}^{d+1}$, $v \in \mathbb{S}^{d}$ and $r>0$ there holds
    \begin{equation}
    \label{eq: tp-radius double sphere equality}
        \left\{ z \in \mathbb{R}^{d+1} \setminus \{x\}| \, \frac{|\langle z-x,v \rangle|}{|z-x|^2}=\frac{1}{2r} \right\}
    =
    (S_{x-rv,r } \cup S_{x+rv,r }) \setminus \{x\}  .
    \end{equation}
    Consequently, if $\Omega \subset \mathbb{R}^{d+1}$ is a bounded Caccioppoli set, then for all $x,y \in \partial ^\ast \Omega$ with $x \neq y$ the tangent-point radius $r_{\TP}(x,y)$ satisfies 
    \begin{equation}
    	\label{eq: def tp-radius}
    	\frac{| \langle x-y,n(x) \rangle |}{|x-y|^2}= \frac{1}{2 r_{\TP}(x,y)} .
    \end{equation}
\end{lem}

\begin{proof}
	Let $x \in \mathbb{R}^{d+1}$, $v \in \mathbb{S}^{d}$, $r>0$ and $z \in \mathbb{R}^{d+1} \setminus \{x\}$ and set $a:=\textnormal{sgn}(\langle z-x,v\rangle)$. If $a=0$, then $z$ lies in neither of the sets in \eqref{eq: tp-radius double sphere equality}. Otherwise,
	\begin{align*}
		&\frac{|\langle z-x,v \rangle|}{|z-x|^2}=\frac{1}{2r}
		\\
		\Leftrightarrow \quad
		&|z-x|^2= 2r | \langle z-x,v\rangle |
		\\
		\Leftrightarrow \quad
		&|z-x|^2-2r|\langle z-x,v \rangle | + r^2 |v|^2 = r^2
		\\
		\Leftrightarrow \quad 
		&|z-x-arv|^2=r^2
		\\
		\Leftrightarrow \quad
		&z \in S_{x+arv,r}
		\\
		\Leftrightarrow \quad
		&z \in S_{x-rv,r } \cup S_{x+rv,r } \, .
	\end{align*}
	This proves \eqref{eq: tp-radius double sphere equality}.
	
	Now let $\Omega \subset \mathbb{R}^{d+1}$ be a bounded Caccioppoli set and $x,y \in \partial ^\ast \Omega$ with $x \neq y$. If $\langle y-x, n(x) \rangle =0$, then by definition $r_{\TP}(x,y) = \infty$ and thus it satisfies \eqref{eq: def tp-radius}. Otherwise, by definition $r_{\TP}(x,y)$ is the radius of the unique sphere $S$ through $x$ and $y$, which is tangential to $\partial ^\ast \Omega$ at $x$. We have that both $S_{x-r n(x),r}$ and $S_{x+r n(x),r}$ contain $x$ and by definition of $n(x)$ are also tangential to $\partial ^\ast \Omega$ at $x$ for every $r>0$. Since $\langle y-x, n(x) \rangle \neq 0$, there exists a unique $r_0>0$ such that $y \in S_{x-r_0 n(x),r_0} \cup S_{x+r_0 n(x),r_0}$. Hence, $S$ coincides with either $S_{x-r_0 n(x),r_0}$ or $S_{x+r_0 n(x),r_0}$ and so $r_{\TP}(x,y)=r_0$. From \eqref{eq: tp-radius double sphere equality} we infer that $r_{\TP}(x,y)$ satisfies \eqref{eq: def tp-radius}.
\end{proof}

In order to prove the rigidity of spheres in \eqref{eq: main theorem}, we first establish another rigidity statement regarding boundaries of convex sets with constant tangent-point radius. 

\begin{lem}
\label{lem: sphere iff constant tangent point radius rectifiable}
    Let $\Omega \subset \mathbb{R}^{d+1}$ be a bounded convex set, $x \in \partial ^\ast \Omega$, $r>0$ and $S:=r \mathbb{S}^{d}$. Then $\partial ^\ast \Omega= \partial  \Omega=S$ (up to translations), if and only if 
    \[
    2 \, \frac{|\langle y-x,n(x) \rangle |}{|x-y|^2}=\frac{1}{r} \quad \textnormal{for } \mathcal{H}^{d}-\textnormal{a.e. } y \in \partial ^\ast \Omega \quad \textnormal{and } \mathcal{H}^{d}(\partial ^\ast \Omega)=\mathcal{H}^{d}(S).
    \]
\end{lem}

\begin{proof}
    If $\partial \Omega$ is a translation of $S$, then in particular $\mathcal{H}^{d}(\partial ^\ast \Omega)=\mathcal{H}^{d}(S)$. Moreover, $\partial \Omega$ coincides for all $x,y \in \partial ^\ast \Omega$ with $x \neq y$ with the unique sphere through $x$ and $y$, which is tangent to $\partial ^\ast \Omega$ at $x$. Therefore, by \Cref{lem: formula for tp-radius} we have
    \[
    2 \, \frac{|\langle y-x,n(x) \rangle |}{|x-y|^2}= \frac{1}{r_{\TP}(x,y)}= \frac{1}{r}  \quad \forall x,y \in \partial ^\ast \Omega , x \neq y \, ,
    \]
    which proves "$\Rightarrow$".
    
    It remains to show "$\Leftarrow$".
    W.l.o.g. let $x=0$ and $n(x)=e_1$. Further define $S_1:= S_{-re_{1},r }$ and $B_1:=B_{-re_{1},r}$.
    Suppose
    \[
     2 \, \frac{|\langle y,e_1 \rangle |}{|y|^2}=\frac{1}{r} \quad \textnormal{for } \mathcal{H}^{d}-\textnormal{a.e. } y \in \partial ^\ast \Omega \quad \textnormal{and } \mathcal{H}^{d}(\partial ^\ast \Omega)= \mathcal{H}^{d}(S) .
    \]
    
     From the first part of \Cref{lem: formula for tp-radius} and the formula in \Cref{prop: Supporting Hyperplane Theorem strong} (i) applied to $x=0$, we infer that there exists $N \subset \mathbb{R}^{d+1}$ with $\mathcal{H}^{d}(N)=0$, such that 
    \[
    \partial ^\ast \Omega \setminus N \subseteq (S_{re_{1},r } \cup S_{-re_{1},r }) \cap ((-\infty ,0] \times \mathbb{R}^d) = S_1.
    \]
  Combined with $\mathcal{H}^{d}(\partial ^\ast \Omega)=\mathcal{H}^{d}(S)= \mathcal{H}^{d}(S_1)$, this yields $\mathcal{H}^d(S_1 \setminus \partial ^\ast \Omega)=0$. So, $\partial ^\ast \Omega \cap S_1$ as well as $\partial \Omega \cap S_1$ are dense in $S_1$. Closedness of $\partial \Omega$ implies that $S_1 \subseteq \partial \Omega$. Thus, the convexity of $\Omega$ further implies that $B_1 \subseteq \Omega$, which by \cite[Theorem~15.5]{FinitePerimeterSets} forces $n_\Omega(x)=n_{B_1}(x)$ for every $x  \in S_1 \cap \partial ^\ast \Omega$. On the other hand, by applying the formula in \Cref{prop: Supporting Hyperplane Theorem strong} (i) to every $x \in S_1 \cap \partial ^\ast \Omega$ using the density of $\partial  ^\ast \Omega \cap S_1$ in $S_1$, we obtain
  \[
  B_1 \subseteq \Omega \subseteq \bigcap_{x \in \partial ^\ast \Omega \cap S_1} 
  \{ y \in \mathbb{R}^{d+1} \, | \, \langle y-x, n_{B_1}(x) \rangle \leq 0\} = 
  \overline{B}_1.
  \]
  Therefore, the symmetric difference between $\Omega$ and $ B_1$ has Lebesgue measure zero, so by \cite[Remark~15.2]{FinitePerimeterSets} we conclude $\partial ^\ast \Omega= S_1= S-re_{1}$ and subsequently $\partial \Omega=S-re_{1}$.
\end{proof}

\begin{cor}
\label{cor: sphere constant integrand}
    Let $d \in \mathbb{N}, q >0$ and let $S$ denote a $d$-dimensional sphere with radius $r$. Then for the integrand of $\TP^{(2q,q)}(S)$ we have
    \[
    \left(\frac{|\langle y-x,n(x) \rangle |}{|x-y|^2} \right)^q=\left(\frac{1}{2r}\right)^q \quad \forall x,y \in S , x \neq y .
    \]
    Let $r_d$ be defined as in \eqref{eq: def r_d}. Then in particular, there holds 
    \begin{equation}
    \label{eq: main theorem lower bound via r}
    \frac{\TP _{q} (S)}{\mathcal{H}^d(S)^{\frac{2d-q}{d}}}= \frac{1}{r_d^q} .
    \end{equation}
\end{cor}

\begin{proof}
    As an immediate consequence of \Cref{lem: sphere iff constant tangent point radius rectifiable}, we obtain that the integrand of $\TP^{(2q,q)}(S)$ is $\mathcal{H}^{2d}$-a.e. equal to $(2r)^{-q}$, which due to continuity reasons can be improved to equality on $S \times S \setminus \{x=y\}$.
    To prove \eqref{eq: main theorem lower bound via r}, we first use the scale invariance of both sides of the equation, allowing us to assume w.l.o.g. that $\mathcal{H}^d(S)=1$. Then by definition $r=r_d$. So 
    \begin{align*}
        \frac{\TP _{q} (S)}{\mathcal{H}^d(S)^{\frac{2d-q}{d}}}
     =
     \TP _{q} (S)
     =
     2^q \, \TP^{(2q,q)} (S)
     =
     2^q \int_S \int_{S} (2r_d)^{-q} \, d \mathcal{H}^d(y) \, d \mathcal{H}^d(x)
     =
     \frac{1}{r_d^q} .
    \end{align*}
\end{proof}

\section{Sharp isoperimetric inequality \\ for classical tangent-point energies \texorpdfstring{{\boldmath$\TP_q$}}{TP_q}}
\label{sec: Proof of the main theorem}

In \cite[Lemma~1]{bushling2024singular} Bushling computed the tangential Jacobian of the Gauss map, defined in \eqref{eq: def: Gauss map}, on $d$-rectifiable sets in $\mathbb{R}^{d+1}$, which reads as follows in our notation.

\begin{lem}
\label{lem: relation between phi and integrand for sets}
   Let $\Upsilon \subset \mathbb{R}^{d+1}$ be a $d$-rectifiable set, $n: \Upsilon \rightarrow \mathbb{S}^{d}$ a measurable unit-normal vector field, $x \in \Upsilon$ a point and $\phi_x$ the Gauss map of $\Upsilon$ and $x$. Then the a.e.-defined Jacobian determinant $|J\phi_x|:\Upsilon \setminus \{x\} \rightarrow \mathbb{R}$ of $\phi_x$ is given by
   \begin{equation}
   \label{eq: Gauss map Jacobian}
   |J \phi_x (y)|= \frac{|\langle x-y,n(y) \rangle |}{|x-y|^{d+1}}  .
   \end{equation}
\end{lem}

Since the formula \eqref{eq: Gauss map Jacobian} is applicable for the reduced boundary of finite-perimeter sets, using the resemblance between \eqref{eq: Gauss map Jacobian} and the integrand of $\TP^{(p,q)}$, we could rewrite the tangent-point energy as
\begin{equation}
\label{eq: TP^{(p,q)} via Gauss map}
\TP^{(p,q)}(\Sigma)= 
\int_{\Sigma^\ast} \int_{\Sigma ^\ast}
|J \phi_x (y) |^q \, \cdot \,
|x-y|^{(d+1)q-p} 
\, d \mathcal{H}^d(y) \, d \mathcal{H}^d(x)  .
\end{equation}
For a quick motivation of the utility of the Gauss map, we may consider the simpler expression $\int_{\Sigma ^\ast} |J \phi_x (y) | \, d \mathcal{H}^{d}(y)$.
Via an application of the coarea formula (see \Cref{prop: coarea formula}) to $\phi_x$ and \Cref{lem: surjectivity of Gauss map}, we obtain 
\[
\int_{\Sigma ^\ast} |J \phi_x (y) | \, d \mathcal{H}^{d}(y) = \int_{\mathbb{S}^{d}} \bigg( \sum_{y \in \phi_x ^{-1}(\omega)} 1 \bigg) \, d \mathcal{H}^{d} (\omega) \geq \mathcal{H}^{d} (\textnormal{Im}(\phi_x)) \geq \frac{1}{2} \mathcal{H}^{d} (\mathbb{S}^{d}) .
\]
In fact, for boundaries $\Sigma$ of convex sets equality holds for every $x \in \Sigma ^\ast$ in both inequalities above. 

However, the right-hand side of \eqref{eq: TP^{(p,q)} via Gauss map} does not have such a convenient form for establishing sharp lower bounds. This is because in general the function $h:=|J \phi_x |^{q-1} \, \cdot \, |x-\cdot|^{(d+1)q-p}$ does not depend on $\phi_x$ in a way that applying \Cref{prop: coarea formula} to $g=h$ and $f=\phi_x$ would simplify \eqref{eq: TP^{(p,q)} via Gauss map}.

Nevertheless, using the constancy of the integrand of $\TP_q=2^q\TP^{(2q,q)}$ on spheres (see \Cref{cor: sphere constant integrand}), in \Cref{subsec: Sharp Minorants of TP as special cases of F_G} we will reduce the minimization problem of $\TP_q$ to one for a functional $F_d$, which can be bounded from below in the fashion of the short motivation.
In fact, $F_d$ belongs to a family of functionals $F_G$, which depend on a measurable function $G:[0,1] \rightarrow [0, \infty] $ and to which the techniques of \cite{bushling2024singular} extend quite naturally. More concretely, $F_d$ corresponds to the choice $G(t):=t^{d-1}$, where $F_G$ is given by
\begin{equation}
\label{eq: F_G definition}
F_{G}(\Sigma)=\int_{\Sigma ^\ast} \int_{\Sigma ^\ast}
G\left(\frac{| \langle x-y,n(x) \rangle |}{|x-y|} \right) \,
\frac{| \langle x-y,n(y) \rangle |}{|x-y|^{d+1}} 
\, d \mathcal{H}^d(y) \, d \mathcal{H}^d(x)  .
\end{equation}
Equivalently, we may express the integrand of $F_G$ via the Gauss map, which yields
\begin{equation}
\label{eq: F_G via phi}
F_{G}(\Sigma)=\int_{\Sigma ^\ast} \int _{\Sigma ^\ast} G (|\langle \phi_x(y), n(x) \rangle |) \, |J \phi_x (y)| \, d \mathcal{H}^{d} (y) \, d \mathcal{H}^d(x) .
\end{equation}

\subsection{A non-local integral identity for the perimeter}
\label{subsec: Global Minimizers of F_G}

This subsection is devoted to showing that $F_G/ \rm{Per}$ is uniquely minimized by boundaries $\Sigma$ of convex sets $\Omega \subset \mathbb{R}^{d+1}$. To this end, setting $a(x,y):=\textnormal{sgn} \langle x-y, n(x) \rangle $ we consider a \textit{signed version} and sharp lower bound of $F_G$ given by 
\begin{equation}
\label{eq: tilde F_G definition}
\tilde{F}_G(\Sigma)
		:=
		\int_{\Sigma ^\ast} \int_{\Sigma ^\ast}
		G\left(\frac{| \langle x-y,n(x) \rangle |}{|x-y|} \right) \, a(x,y) \, 
		\frac{\langle y-x,n(y) \rangle }{|x-y|^{d+1}} 
		\, d \mathcal{H}^{d}(y) \, d \mathcal{H}^{d}(x)  .
\end{equation}
Mimicking the proof of \cite[Theorem, p.~3]{bushling2024singular}, we show that (if well defined) $\tilde{F}_G$ is actually a constant multiple of the perimeter on the domain of bounded Caccioppoli sets.

\begin{thm}
\label{thm: bushling theorem general}
    Let $d \in \mathbb{N}$, $G:[0,1] \rightarrow [0, \infty]$ be measurable and bounded and $\Sigma$ be the boundary of a bounded Caccioppoli set $\Omega \subset \mathbb{R}^{d+1}$. Assume that $F_G(\Sigma)< \infty$. Then for $\mathcal{H}^d$-a.e. $x \in \Sigma^\ast=\partial ^\ast \Omega$ the identity 
    \begin{equation}
    	\label{eq: bushling general}
        \int_{\Sigma ^\ast }
        G\left(\frac{| \langle x-y,n(x) \rangle |}{|x-y|} \right) \, \textnormal{sgn}  \langle x-y, n(x) \rangle \, 
        \frac{\langle y-x,n(y) \rangle }{|x-y|^{d+1}} 
        \, d \mathcal{H}^d(y)= b_{d,G}
    \end{equation}
    holds with $b_{d,G}:= \frac{1}{2} \, 
    \int_{\mathbb{S}^d}
    G\left(| \langle \omega ,e_1 \rangle |\right) \, 
    d\mathcal{H}^d ( \omega)$.
    Consequently, $\tilde{F}_G(\Sigma)=b_{d,G} \, \mathcal{H}^{d}(\Sigma ^\ast)= b_{d,G} \, \textnormal{Per}(\Omega)$.
\end{thm}

\begin{proof}
    According to \cite[Lemma~2]{bushling2024singular}, $\Sigma ^\ast$ satisfies the OCC. Fix a point $x \in \Sigma^\ast$ at which the OCC is satisfied and for which 
    \begin{equation}
    \label{eq: bushling help}
    \int_{\Sigma ^\ast }
    G\left(\frac{| \langle x-y,n(x) \rangle |}{|x-y|} \right) \,
    \frac{|\langle y-x,n(y) \rangle |}{|x-y|^{d+1}} 
    \, d \mathcal{H}^d(y) < \infty 
    \end{equation}
    holds, which is the case for $\mathcal{H}^{d}$-a.e. $x \in \Sigma ^\ast$ by the assumption $F_G(\Sigma)< \infty$ and Tonelli's theorem. Let $L_{x,\omega}:=x + \textnormal{span } \omega$. Then by applying the coarea formula (\Cref{prop: coarea formula}) to $U=\Sigma^\ast \setminus \{x\}$, $V=\mathbb{S}^{d}$, $f= \phi_x$ and
    \[
    g: \Sigma^\ast \setminus \{x\} \rightarrow \R , \quad 
    g(y):=G\left(| \langle \phi_x (y) ,n(x) \rangle | \right) \, \textnormal{sgn }  \langle -\phi_x(y) , n(x) \rangle \, \textnormal{sgn} \langle \phi_x(y), n(y) \rangle ,
    \]
    we get
    \begin{align*}
    	&\int_{\Sigma ^\ast}
    	G\left(\frac{| \langle x-y,n(x) \rangle |}{|x-y|} \right) \, \textnormal{sgn }  \langle x-y, n(x) \rangle \, 
    	\frac{\langle y-x,n(y) \rangle }{|x-y|^{d+1}} 
    	\, d \mathcal{H}^d(y)
    	\\
    	=&
        \int_{\Sigma ^\ast} G\left(| \langle \phi_x (y) ,n(x) \rangle | \right) \, \textnormal{sgn }  \langle -\phi_x(y) , n(x) \rangle \, \textnormal{sgn} \langle \phi_x(y), n(y) \rangle \, |J \phi_x(y)| 
        \, d \mathcal{H}^d(y)
        \\
        =&
        \int_{\mathbb{S}^d} \sum_{\phi_x^{-1}(\omega)}
        G\left(| \langle \omega ,n(x) \rangle |\right) \, 
        \, \textnormal{sgn } \langle - \omega, n(x) \rangle \, 
        \textnormal{sgn } \langle  \omega , n(y) \rangle 
        \, d\mathcal{H}^d ( \omega)
        \\
        =&
        \int_{\mathbb{S}^d} \sum_{\overset{\lambda >0}{x + \lambda \omega \in \Sigma ^\ast }}
        G\left(| \langle \omega ,n(x) \rangle |\right) \, 
        \textnormal{sgn } \langle - \omega , n(x) \rangle \,
        \textnormal{sgn } \langle  \omega , n(x + \lambda \omega ) \rangle \, 
        d\mathcal{H}^d ( \omega) 
        \\
        =&
        \frac{1}{2} \, 
        \int_{\mathbb{S}^d} \sum_{\overset{\lambda >0}{x + \lambda \omega \in \Sigma ^\ast }}
        G\left(| \langle \omega ,n(x) \rangle |\right) \, 
        \textnormal{sgn } \langle - \omega , n(x) \rangle \,
        \textnormal{sgn } \langle  \omega , n(x + \lambda \omega ) \rangle \, 
        d\mathcal{H}^d ( \omega) 
        \\
        +&
        \frac{1}{2} \, 
        \int_{\mathbb{S}^d} \sum_{\overset{\lambda >0}{x + \lambda (-\omega) \in \Sigma ^\ast }}
        G\left(| \langle -\omega ,n(x) \rangle |\right) \, 
        \textnormal{sgn } \langle \omega , n(x) \rangle \,
        \textnormal{sgn } \langle  - \omega , n(x + \lambda (-\omega) ) \rangle \, 
        d\mathcal{H}^d ( \omega) 
        \\
        =&
        \frac{1}{2} \,
        \int_{\mathbb{S}^d}
        \sum_{y \in \Sigma ^\ast \cap (L_{x,\omega} \setminus \{x\})}
        G\left(| \langle \omega ,n(x) \rangle |\right) \, 
        \textnormal{sgn } \langle - \omega , n(x) \rangle \,
        \textnormal{sgn } \langle  \omega , n(y) \rangle \, 
        d\mathcal{H}^d ( \omega) .
    \end{align*}
    By the orientation cancellation condition,
    \begin{align*}
        \sum_{y \in \Sigma ^\ast \cap (L_{x,\omega} \setminus \{x\})}
        \textnormal{sgn }\langle \omega  , n(y) \rangle
        =&
        \left( \sum_{y \in \Sigma ^\ast \cap L_{x,\omega} } \textnormal{sgn }\langle \omega  , n(y) \rangle \right) 
        \, - \, 
        \textnormal{sgn } \langle \omega , n(x) \rangle 
        \\
        =&
        \textnormal{sgn } \langle - \omega , n(x) \rangle \, ,
    \end{align*}
     so we conclude that
     \begin{align*}
          &\int_{\Sigma ^\ast }
        G\left(\frac{| \langle x-y,n(x) \rangle |}{|x-y|} \right) \, \textnormal{sgn }  \langle x-y, n(x) \rangle \, 
        \frac{\langle y-x,n(y) \rangle }{|x-y|^{d+1}} 
        \, d \mathcal{H}^d(y)
        \\
        =&
        \frac{1}{2} \, 
        \int_{\mathbb{S}^d}
        G\left(| \langle \omega ,n(x) \rangle |\right) \, 
        \textnormal{sgn } \langle - \omega , n(x) \rangle \,
        \textnormal{sgn } \langle  - \omega , n(x) \rangle \, 
        d\mathcal{H}^d ( \omega)
        \\
        =&
        \frac{1}{2} \, 
        \int_{\mathbb{S}^d}
        G\left(| \langle \omega ,n(x) \rangle |\right) \, 
        d\mathcal{H}^d ( \omega) .
     \end{align*}
     This last integral is invariant under rotations and therefore
     \begin{align*}
     &\int_{\Sigma ^\ast }
        G\left(\frac{| \langle x-y,n(x) \rangle |}{|x-y|} \right) \, \textnormal{sgn }  \langle x-y, n(x) \rangle \, 
        \frac{\langle y-x,n(y) \rangle }{|x-y|^{d+1}} 
        \, d \mathcal{H}^d(y)
        \\
        =&
        \frac{1}{2} \, 
        \int_{\mathbb{S}^d}
        G\left(| \langle \omega ,e_1 \rangle |\right) \, 
        d\mathcal{H}^d ( \omega)=b_{d,G} .
     \end{align*}
    By integrating with respect to $d \mathcal{H}^d (x)$, from the definition of $\tilde{F}_G$ \eqref{eq: tilde F_G definition} we immediately conclude $\tilde{F}_G= b_{d,G} \mathcal{H}^{d}(\Sigma ^\ast)$.
\end{proof}

\begin{rmk}
	As in \cite{bushling2024singular}, the first part of \Cref{thm: bushling theorem general} holds more generally for pairs $(\Upsilon, n)$ satisfying the OCC and the integrability hypothesis $F_G(\Upsilon)<\infty$.
\end{rmk}

As a consequence of \Cref{thm: bushling theorem general}, we obtain a generalization of \cite[Theorem, p.~1]{steinerberger2022inequality} and \cite[Corollary, p.~4]{bushling2024singular}.
\begin{cor}
\label{cor: lower bound for F_G for sets}
    Let $d \in \mathbb{N}$, $G:[0,1] \rightarrow [0, \infty]$ be measurable and bounded, satisfying $G(t)>0$ for every $t \in (0,1]$. Furthermore, let $b_{d,G}$ be defined as in \Cref{thm: bushling theorem general} and $\Sigma$ be the boundary of a bounded Caccioppoli set $\Omega \subset \mathbb{R}^{d+1}$ with $\textnormal{Per}(\Omega)>0$. Then
    \[
    \frac{F_{G}(\Sigma)}{\mathcal{H}^d (\Sigma ^\ast)} \geq \frac{F_G(S)}{\mathcal{H}^d(S)}=b_{d,G} ,
    \]
    where $S$ denotes a $d$-dimensional sphere. Here, equality holds if and only if  $\Sigma ^\ast$ is the reduced boundary of a convex set.
\end{cor}
\begin{proof}
    Let $G, \Sigma$ be given as above. W.l.o.g. $\mathcal{H}^d(\Sigma ^\ast) = \mathcal{H}^d(S)=1$. If $F_G(\Sigma)=\infty$, the bound is immediate, so assume $F_G(\Sigma)<\infty$. In particular, this is the case if $\Omega$ is convex, as the following argument shows. By \eqref{eq: F_G via phi}, boundedness of $G$ and \Cref{prop: coarea formula} we get 
    \[
    F_G(\Sigma)\leq M \int_{\Sigma ^\ast} \int_{\Sigma ^\ast}
    |J\phi_x(y)|
    \, d \mathcal{H}^{d}(y) d \mathcal{H}^{d} (x)
    =
     M \int_{\Sigma ^\ast} \int_{\mathbb{S}^{d}}
    \mathcal{H}^{0} (\phi_x ^{-1}(\omega))
    \, d \mathcal{H}^{d}(\omega) d \mathcal{H}^{d} (x) 
    \] 
    with $M:=\sup_{t \in [0,1]} G(t)< \infty$. If $\Omega$ is convex, then by \Cref{lem: surjectivity of Gauss map} we have $\mathcal{H}^{0} (\phi_x ^{-1}(\omega))\leq 1$ for $\mathcal{H}^{d}$-a.e. $x \in \Sigma ^\ast$ and $\omega \in \mathbb{S}^{d}$, so $F_G(\Sigma)\leq M \mathcal{H}^{d}(\mathbb{S}^{d})< \infty$.
    
    Denote by $f_G$ the integrand of $F_G(\Sigma)$ and by $\tilde{f}_G$ the integrand of $\tilde{F}_G(\Sigma)$ (see \eqref{eq: tilde F_G definition}). Since $f_G \geq \tilde{f}_G$, we have by \Cref{thm: bushling theorem general} $F_G(\Sigma) \geq \tilde{F}_G(\Sigma)=b_{d,G}=F_G(S)$ with equality if and only if  $f_G(x,y)=\tilde{f}_G(x,y)$ for  $\mathcal{H}^{2d}$-a.e. $(x,y) \in \Sigma ^\ast \times \Sigma ^\ast$. 
     
     Define $A:=\{(x,y) \in  \Sigma ^\ast \times \Sigma ^\ast \, | \, \langle x-y, n(x) \rangle =0 \}$. 
    
     If $(x,y) \notin A$, then $f_G(x,y)= |\tilde{f}_G(x,y)|$, so
     \[
     f_G(x,y)= \tilde{f}_G(x,y) 
     \quad \Leftrightarrow \quad 
     \tilde{f}_G(x,y) \geq 0 
     \quad \Leftrightarrow \quad 
     \langle x-y , n (y) \rangle \langle y-x, n (x) \rangle \geq 0 \, ,
     \]
     since by assumption $G(t)>0$ for every $t \in (0,1]$.
     
     On the other hand, $f_G(x,y)=\tilde{f}_G(x,y)=0$ for $\mathcal{H}^{2d}$-a.e. $(x,y) \in A$, because $|\tilde{f}_G(x,y)| \leq |f_G(x,y)| \leq M \, |J \phi_x(y)|$ and by \Cref{prop: coarea formula}
    \[
    \int_A |J \phi_x(y)|
     \, d \mathcal{H}^d(y) \, d \mathcal{H}^d(x)
     =
    \int_{\Sigma^\ast} \int_{\mathbb{S}^{d}} \mathcal{H}^{0}(\phi_x^{-1}(\omega)) \, \chi_{\{ \langle \cdot , n(x) \rangle =0 \}} (\omega)
     \, d \mathcal{H}^d(\omega) \, d \mathcal{H}^d(x) =0 \, ,
    \] 
    as the $\mathcal{H}^{d}$-measure of the equator $\{ \omega \in \mathbb{S}^{d} \, | \, \langle \omega , n(x) \rangle =0 \}$ is zero.
    
    Therefore, equality of $F_G(\Sigma)$ and $\tilde{F}_G(\Sigma)$ is equivalent to
    \[
    \langle x-y , n (y) \rangle \langle y-x, n (x) \rangle \geq 0 \quad \textnormal{for } \mathcal{H}^{2d}-\textnormal{a.e. } (x,y) \in \Sigma ^\ast \times \Sigma ^\ast \, ,
    \]
    which by \Cref{lem: convex set nonlocal characterization} is equivalent to $\Sigma ^\ast$ being the reduced boundary of a convex set. 
\end{proof}

As the function $G_\varrho:t \mapsto t^{\varrho-1}$ satisfies the assumptions from \Cref{cor: lower bound for F_G for sets} for all $\varrho\geq 1$, for $F_\varrho:=F_{G_\varrho}$, that is
\begin{equation}
\label{eq: F_q definition}
	F_\varrho(\Sigma)= 
	\int_{\Sigma ^\ast }\int_{\Sigma ^\ast}
	\left(\frac{| \langle x-y,n(x) \rangle |}{|x-y|} \right)^{\varrho-1} \,
	\frac{| \langle x-y,n(y) \rangle |}{|x-y|^{d+1}} 
	\, d \mathcal{H}^d(y) \, d \mathcal{H}^d(x)  \, ,
\end{equation}
we immediately conclude the following.
\begin{cor}
\label{cor: lower bound on F_q for sets}
    Let $d \in \mathbb{N}$, $\varrho \in [1, \infty)$, $b_{d,\varrho}:=b_{d,G_\varrho}= \frac{1}{2} \, 
    \int_{\mathbb{S}^d}
    | \langle \omega ,e_1 \rangle |^{\varrho -1} \, 
    d\mathcal{H}^d ( \omega)$  and $\Sigma$ be the boundary of a bounded Caccioppoli set $\Omega \subset \mathbb{R}^{d+1}$ with $\textnormal{Per}(\Omega)>0$. Then 
   \[
    \frac{F_\varrho(\Sigma )}{\mathcal{H}^{d}(\Sigma ^\ast)} \geq  \, \frac{F_\varrho(S)}{\mathcal{H}^d(S)}= b_{d, \varrho},
    \]
    where $S$ denotes a $d$-dimensional sphere. Here, equality holds if and only if $\Sigma ^\ast$ is the reduced boundary of a convex set.
\end{cor}

\subsection{Sharp minorants of the tangent-point energies}
\label{subsec: Sharp Minorants of TP as special cases of F_G}

We will now translate the results from the previous subsection to $\TP^{(q+d,q)}$ and subsequently to $\TP_q$, finishing the proof of \Cref{thm: main 1}.

By a symmetrization argument, consisting of a combination of Tonelli's theorem and Young's inequality, applied to the integrand of $F_\varrho$, we obtain the same sharp lower bound for $\TP^{(q+d,q)}$ as for $F_\varrho$ with the choice $\varrho=q$.

\begin{thm}
\label{thm: lower bound on TP_q+d for sets}
    Let $d \in \mathbb{N}$, $q \in [1, \infty)$, $b_{d,q}= \frac{1}{2} \, 
    \int_{\mathbb{S}^d}
    | \langle \omega ,e_1 \rangle |^{q -1} \, 
    d\mathcal{H}^d ( \omega)$ and $\Sigma$ be the boundary of a bounded Caccioppoli set $\Omega \subset \mathbb{R}^{d+1}$ with $\textnormal{Per}(\Omega)>0$. Then 
    \[
    \frac{\TP^{(q+d,q)}(\Sigma)}{\mathcal{H}^{d}(\Sigma ^\ast)} \geq 
    \frac{\TP^{(q+d,q)}(S)}{\mathcal{H}^d(S)}= b_{d,q},
    \]
    where $S$ denotes a $d$-dimensional sphere. If equality holds, then $\Sigma ^\ast$ is the reduced boundary of a convex set.
\end{thm}

\begin{proof}
     W.l.o.g. $\mathcal{H}^d(\Sigma ^\ast)=\mathcal{H}^d(S)=1$. Notice that $\TP^{(1+d,1)}=F_1$. Otherwise, for $q>1$ we have by Tonelli's theorem, Young's inequality and \Cref{cor: lower bound on F_q for sets}
    \begin{align*}
        &\TP^{(q+d,q)}(\Sigma)
        \\
        =&
        \frac{q-1}{q} \,
        \int_{\Sigma ^\ast} \int_{\Sigma ^\ast} \frac{| \langle x-y,n(x) \rangle |^q}{|x-y|^{q+d}} \, d \mathcal{H}^d(y) \, d \mathcal{H}^d(x)
        \\
        +&
        \frac{1}{q} \, 
        \int_{\Sigma ^\ast} \int_{\Sigma ^\ast} \frac{| \langle x-y,n(y) \rangle |^q}{|x-y|^{q+d}} \, d \mathcal{H}^d(y) \, d \mathcal{H}^d(x)
        \\
        =&
         \int_{\Sigma ^\ast} \int_{\Sigma ^\ast}
        \frac{\frac{q-1}{q}| \langle x-y, n(x) \rangle |^q \, + \, \frac{1}{q} | \langle x-y, n(y) \rangle |^q}{|x-y|^{q+d}} 
        \, d \mathcal{H}^{d}(y) \,  d \mathcal{H}^{d}(x)
        \\
        \geq&
        \int_{\Sigma ^\ast} \int_{\Sigma ^\ast}
        \frac{| \langle x-y,n(x) \rangle |^{q-1} \, | \langle x-y,n(y) \rangle |}{|x-y|^{q+d}}
        \, d \mathcal{H}^d(y) \, d \mathcal{H}^d(x)
       	\\
       	=&
       	F_q(\Sigma) \, 
    \end{align*}
    with equality for spheres since 
   \[
   | \langle x-y,n(x) \rangle |= | \langle x-y,n(y) \rangle | \quad \forall x,y \in S .
   \]
   Furthermore by \Cref{cor: lower bound on F_q for sets} for all $q \geq 1$ there holds $F_q(\Sigma) \geq F_q(S)= \TP^{(q+d,q)}(S)=b_{d,q}$ with equality if and only if $\Sigma ^\ast$ is the reduced boundary of a convex set, which concludes the proof.
\end{proof}

On the other hand, from the previous section we can also deduce an upper bound for $\TP^{(q+d,q)}$ on the domain of boundaries $\Sigma$ of bounded convex sets $\Omega$ with fixed perimeter, as shown in the following corollary. In particular, the energy does not blow up for sequences $(\Sigma_k)_k$ of such boundaries for which $\Omega_k$ converges in Hausdorff distance to a lower-dimensional convex set. While \Cref{cor: upper bound on TP} is not used in the rest of the paper, it helps to distinguish the characteristics of the lower limit case from those of the subcritical case.

\begin{cor}
\label{cor: upper bound on TP}
    Let $d \in \mathbb{N}$, $q \in [1, \infty)$ and $\Sigma$ be the boundary of a bounded convex set $\Omega \subset \mathbb{R}^{d+1}$ with $\textnormal{Per}(\Omega)=\mathcal{H}^{d}(\Sigma ^\ast)=\mathcal{H}^{d}(\Sigma)>0$. Then 
    \[
    \frac{\TP^{(q+d,q)}(\Sigma)}{\mathcal{H}^d(\Sigma)} \leq \frac{\TP^{(1+d,1)}(S)}{\mathcal{H}^d(S)}=  \frac{1}{2} \mathcal{H}^d(\mathbb{S}^d) \, ,
    \]
    where $S$ denotes a $d$-dimensional sphere and $\mathbb{S}^d$ denotes the $d$-dimensional unit sphere.
\end{cor}

\begin{proof}
    W.l.o.g. $\mathcal{H}^d(\Sigma)=\mathcal{H}^d(S)=1$.
    We have 
    \begin{align*}
    \TP^{(q+d,q)}(\Sigma)=& 
    \int_{\Sigma ^\ast} \int_{\Sigma ^\ast} \frac{| \langle x-y,n(y) \rangle |}{|x-y|^{1+d}} \, 
    \left(\frac{| \langle x-y,n(y) \rangle |}{|x-y|}\right)^{q-1}
    \, d \mathcal{H}^d(y) \, d \mathcal{H}^d(x)
    \\
    \leq&
    \int_{\Sigma ^\ast} \int_{\Sigma ^\ast} \frac{| \langle x-y,n(y) \rangle |}{|x-y|^{1+d}} 
    \, d \mathcal{H}^d(y) \, d \mathcal{H}^d(x)
    = \TP^{(1+d,1)}(\Sigma)  .
    \end{align*}
    On the other hand, $\TP^{(1+d,1)}$ coincides with $F_1$, given by \eqref{eq: F_q definition}, so by \Cref{cor: lower bound on F_q for sets} we get 
    \[
    	\TP^{(1+d,1)}(\Sigma)=\TP^{(1+d,1)}(S)=b_{d,1}= \frac{1}{2} \mathcal{H}^d(\mathbb{S}^d) .
    \]
\end{proof}

As a further consequence, using Jensen's inequality we conclude the proof of \Cref{thm: main 1}.

\begin{proof}[Proof of \Cref{thm: main 1}]

     W.l.o.g. $\mathcal{H}^d(\Sigma ^\ast)=\mathcal{H}^d(S)=1$. Then we have by Jensen's inequality, \Cref{thm: lower bound on TP_q+d for sets} and \Cref{cor: sphere constant integrand}
     \begin{equation}
      \label{eq: main theorem long help}
    \begin{aligned}
          \TP_q(\Sigma)=&
          2^q \int_{\Sigma^\ast} \int_{\Sigma^\ast} 
          \left(\frac{| \langle x-y,n_\Omega(y) \rangle |}{|x-y|^{2}}\right)^q \, d \mathcal{H}^d(y) \, d \mathcal{H}^d(x)
          \\
          \geq&
          2^q \left( 
          \int_{\Sigma ^\ast} \int_{\Sigma ^\ast} 
          \left(\frac{| \langle x-y,n_\Omega(y) \rangle |}{|x-y|^{2}}\right)^d \,
          d \mathcal{H}^d(y) \, d \mathcal{H}^d(x)
          \right)^{\frac{q}{d}}
          \\
          =&
          2^q \left(
          \TP^{(d+d,d)}(\Sigma)
          \right)^{\frac{q}{d}}
          \\
          \geq&
          2^q \left(\TP^{(d+d,d)}(S) \right)^{\frac{q}{d}}
          =
          \left(\TP_d(S) \right)^{\frac{q}{d}}
          =
          \frac{1}{r_d^q}
          = \TP_q(S) .
    \end{aligned}
    \end{equation}
    Now assume that $q>d$ and $\TP_q(\Sigma) = r_d^{-q}$. Then both inequalities in \eqref{eq: main theorem long help} are equalities, so in particular $ \TP^{(d+d,d)}(\Sigma)=\TP^{(d+d,d)}(S)$. We then conclude from \Cref{thm: lower bound on TP_q+d for sets} that $\Sigma ^\ast$ is the reduced boundary of a bounded convex set $\tilde{\Omega}$ with $n_{\Omega}=\pm n_{\tilde{\Omega}}$ by equality of the approximate tangent spaces (see \cite[Corollary~16.1]{FinitePerimeterSets}). Furthermore, equality in Jensen's inequality together with $\mathcal{H}^{d}(\Sigma ^\ast)=1$ implies that the integrand of the first integral in \eqref{eq: main theorem long help} is $\mathcal{H}^{2d}$-a.e. equal to $(2 r_d)^{-q}$, so in particular after interchanging the variables there exists a point $x \in \partial ^\ast \tilde{\Omega}$ such that 
    \[
    2 \, \frac{| \langle x-y,n_{\tilde{\Omega}}(x) \rangle |}{|x-y|^{2}}=
    2 \, \frac{| \langle x-y,n_\Omega(x) \rangle |}{|x-y|^{2}}=\frac{1}{r_d} \quad \textnormal{for } \mathcal{H}^{d}-\textnormal{a.e. } y \in \partial ^\ast \tilde{\Omega} .
    \]
    Hence, we conclude by \Cref{lem: sphere iff constant tangent point radius rectifiable} that $\Sigma ^\ast= \partial ^\ast \tilde{\Omega}$ is a sphere. 
\end{proof}

\section{Sharp isoperimetric inequality \\ for generalized tangent-point energies \texorpdfstring{{\boldmath$\TP^{(p,q)}$}}{TP^(p,q)}}
\label{sec: More parameters}

In the following, we want to extend the sharp lower bounds to the generalized tangent-point energies for the range of $p,q$ being as wide as possible. To this end, we first recall from the introduction the definition of the Riesz energy of a Lebesgue measurable set $\Omega \subset \mathbb{R}^{d+1}$, $d \geq 1$, that is
\begin{equation}
	\label{eq: Riesz energy}
	R_\alpha (\Omega):= \int_{\Omega} \int_{\Omega}
	\frac{1}{|x-y|^{\alpha}} \, dx \, dy \, , \quad \alpha \in (0, d+1)  .
\end{equation}
$R_\alpha(\Omega)$ immediately extends continuously to $\alpha=0$ via $R_0(\Omega):=\mathcal{L}^{d+1}(\Omega)^2$. Moreover, for any bounded, measurable $\Omega \subset \mathbb{R}^{d+1}$ there holds
\begin{equation}
	\label{eq: R_m definition}
	R_{d+1}(\Omega):= \mathcal{H}^{d}(\mathcal{S}^{d}) \, \mathcal{L}^{d+1}(\Omega)= 
	\lim\limits_{\alpha \nearrow d+1} (-\alpha+d+1) \, R_\alpha(\Omega) \, ,
\end{equation}
which can be seen as a consequence of \cite[Theorem~2.1]{RieszLimit} with the choice $u = \chi_\Omega$. 

The following immediate consequence of the well-known isoperimetric inequality for the Riesz energy (see \cite[Equation~(1.2)]{RieszIso}) will be essential for us later.
\begin{prop}
	\label{prop: Riesz maximum}
	Let $d \geq 1$, $\alpha \in  [0,d+1]$ and $\Omega$ be Lebesgue measurable with $0<\mathcal{L}^{d+1}(\Omega)<\infty$. Let $B$ denote a $(d+1)$-dimensional ball. Then
	\begin{equation}
		\label{eq: Riesz maximum}
	\frac{R_\alpha (\Omega)}{\mathcal{L}^{d+1}(\Omega)^{\frac{2(d+1)-\alpha}{d+1}}}
	\leq 
	\frac{R_\alpha (B)}{\mathcal{L}^{d+1}(B)^{\frac{2(d+1)-\alpha}{d+1}}}  .
	\end{equation}
\end{prop}

Via the following identity \eqref{eq: Riesz =< TP sym}, the sharp upper bound on the Riesz energy \eqref{eq: Riesz maximum} translates to one for a double boundary integral, which can be viewed as a signed symmetrized tangent-point functional.
\begin{lem}
\label{lem: TP_sym and Riesz relation}
    Let $d \geq 1$, $\alpha \in [0,d+1]$ and $\Omega \subset \mathbb{R}^{d+1}$ be a bounded Caccioppoli set. Then 
    \begin{equation}
    \label{eq: Riesz =< TP sym}
       c_{d,\alpha} \, R_\alpha(\Omega)= 
    \int_{\partial ^\ast \Omega}
        \int_{\partial ^\ast \Omega}
        \frac{\langle x-y, n(x)\rangle \, \langle y-x, n(y) \rangle}{|y-x|^\alpha} 
        \, d \mathcal{H}^{d}(y) \, d \mathcal{H}^{d}(x) \, ,
    \end{equation}
    	where $c_{d,\alpha}=(-\alpha+d+2)(-\alpha+d+1)$ for $\alpha <d+1$ and $c_{d,d+1}=1$.
\end{lem}
\begin{proof}
    Let $d, \alpha, \Omega$ be given as above. By boundedness of $\Omega$ there exists an $R>0$, such that $\Omega \subset B_{R}(0)$. 
	
	Fix $x \in \partial ^\ast \Omega$. As a consequence of \cite[Corollary~15.8]{FinitePerimeterSets} we have $\mathcal{H}^{d}(B_r(x) \cap \partial ^\ast \Omega) \in \mathcal{O}(r^{d})$ and thus 
	\begin{equation}
		\begin{split}
		&\label{eq: TP_ssym inner integral finite}
	\int_{\partial ^\ast \Omega}
	\frac{|\langle x-y, n(x)\rangle \, \langle y-x, n(y) \rangle|}{|y-x|^\alpha} 
	\, d \mathcal{H}^{d}(y) 
	\leq	\int_{\partial ^\ast \Omega} \frac{1}{|y-x|^{\alpha-2}} 
		\, d \mathcal{H}^{d}(y) < \infty  .
		\end{split}
	\end{equation}
    Define the vector field
    \[
    F_x: \mathbb{R}^{d+1} \setminus\{x\} \rightarrow \mathbb{R}^{d+1}, \quad F_x(y)=\langle x-y, n(x) \rangle \frac{y-x}{|y-x|^\alpha} .
    \]
    Then $F_x \in C^1(\mathbb{R}^{d+1} \setminus \{x\}, \mathbb{R}^{d+1})$ with
    \[
    \textnormal{div } F_x(y)= (-\alpha+d+2) \, \frac{\langle x-y, n(x)\rangle}{|x-y|^\alpha}  .
    \]
    Further define for fixed $\varepsilon>0$ 
    \[
    K_\varepsilon: [0, \infty) \rightarrow [\varepsilon, \infty), \quad K_\varepsilon(r)= (r^2+\varepsilon^2)^{\frac{1}{2}} 
    \]
    and the regularized vector field
    \[
    F_{x, \varepsilon}: \mathbb{R}^{d+1} \rightarrow \mathbb{R}^{d+1}, \quad F_{x, \varepsilon}(y)=\langle x-y, n(x) \rangle \frac{y-x}{K_\varepsilon(|y-x|)^\alpha} .
    \]
    Clearly, $F_{x, \varepsilon} \in C^1(\mathbb{R}^{d+1}, \mathbb{R}^{d+1})$, $|F_{x,\varepsilon}| \leq |F_x|$ and $F_{x, \varepsilon} \overset{\varepsilon \to 0}{\rightarrow} F_x$ in $C^1_{\textnormal{loc}}$ on $\mathbb{R}^{d+1} \setminus \{x\}$. Moreover, 
    \[
    |\textnormal{div } F_{x,\varepsilon}(y)|= 
    \left| \langle x-y,n(x) \rangle  
  	\, \left[ \frac{d+2-\alpha }{K_\varepsilon (|x-y|)^\alpha } \, + \, 
  	\frac{\alpha \varepsilon^2}{K_\varepsilon (|x-y|)^{\alpha +2}} \right] \right| 
  	\leq \frac{d+2}{|x-y|^{\alpha -1}}  .
    \]
    Hence, by applying the Gauss-Green formula (\Cref{prop: Gauss-Green}) to $F_{x, \varepsilon}$ (times a suitable smooth cutoff function with compact support), in combination with dominated convergence supplied by the above bound and \eqref{eq: TP_ssym inner integral finite}, we get
    \begin{align*}
        &
        \int_{\partial ^\ast \Omega}
        \frac{\langle x-y, n(x)\rangle \, \langle y-x, n(y) \rangle}{|y-x|^\alpha} 
        \, d \mathcal{H}^{d}(y) 
        \\
        =&
         \int_{\partial ^\ast \Omega}
        \langle F_x (y) , n(y) \rangle 
        \, d \mathcal{H}^{d}(y) 
        \\
        =&
        \lim\limits_{\varepsilon \to 0}
        \int_{\partial ^\ast \Omega}
        \langle F_{x, \varepsilon} (y) , n(y) \rangle 
        \, d \mathcal{H}^{d}(y) 
        \\
        =&
        \lim\limits_{\varepsilon \to 0}
        \int_{\Omega}
        \textnormal{div } F_{x, \varepsilon} (y) 
        \, d y
        \\
        =&
         \int_{\Omega}
        \textnormal{div } F_{x} (y) 
        \, d y
        \\
        =&
        (- \alpha + d +2) \int_\Omega \frac{\langle x-y, n(x)\rangle}{|x-y|^\alpha}
        \, dy .
    \end{align*}
    Notice that for every $R>1$, $\alpha \in [0, d+1]$, and every Caccioppoli set $\Omega \subset B_R(0)$ there holds
    \begin{equation}
    	\label{eq: Riesz TP_sym intermediate bound}
    \frac{|\langle x-y,n(x)\rangle|}{|x-y|^\alpha}
    \leq
    (2R)^{d+1} |x-y|^{-d}  \quad \textnormal{for all } x \in \partial ^\ast \Omega, y \in \Omega ,
	\end{equation}
	and
	\begin{equation}
	\label{eq: Riesz TP_sym intermediate bound 2}
		\begin{split}
			\int_{\partial ^\ast \Omega} \int_\Omega
			|x-y|^{-d}
			\, dy \, d \mathcal{H}^{d}(x) 
			\leq& 
			\int_{\partial ^\ast \Omega} \int_{B_R(0)}
			|x-y|^{-d} 
			\, dy \, d \mathcal{H}^{d}(x) 
			\\
			\leq&
		    \mathcal{H}^{d} (\partial ^\ast \Omega )
			\, \int_{B_{2R} (0)} |y|^{-d} 
			\, dy < \infty  .
		\end{split}
	\end{equation}
    In particular, the right-hand side in \eqref{eq: Riesz =< TP sym} is thus well defined as an iterated integral. Furthermore, we get by the previous calculations and Fubini's theorem
    \begin{equation}
    	\label{eq: Riesz TP_sym intermediate}
    	\begin{split}
    &\int_{\partial ^\ast \Omega}
    \int_{\partial ^\ast \Omega}
    \frac{\langle x-y, n(x)\rangle \, \langle y-x, n(y) \rangle}{|y-x|^\alpha} 
    \, d \mathcal{H}^{d}(y) \, d \mathcal{H}^{d}(x) 
    \\
    =&
    (-\alpha+d+2) \,
    \int_{\partial ^\ast \Omega} \int_{\Omega}
    \frac{\langle x-y, n(x)\rangle}{|x-y|^\alpha} 
    \, dy \, d \mathcal{H}^{d}(x) 
    \\
    =&
     (-\alpha+d+2) \,
    \int_{\Omega} \int_{\partial ^\ast \Omega}
     \frac{\langle x-y, n(x)\rangle}{|x-y|^\alpha} 
    \, d \mathcal{H}^{d}(x) \, dy .
    	\end{split}
    \end{equation}
    Now assume $\alpha <d+1$ and define for fixed $y \in \mathbb{R}^{d+1}$ 
    the vector field
    \[
    G_y: \mathbb{R}^{d+1} \setminus \{y\} \rightarrow \mathbb{R}^{d+1} , \quad G_y(x)=\frac{x-y}{|x-y|^\alpha}.
    \]
    Then $G_{y} \in C^1(\mathbb{R}^{d+1} \setminus \{y\}, \mathbb{R}^{d+1})$ with 
    \[
    \textnormal{div } G_y (x)= (-\alpha+d+1) \, \frac{1}{|x-y|^\alpha}  .
    \]
    Analogously to the construction for $F_x$, define for fixed $\varepsilon>0$ the regularized vector field
    \[
    G_{y,\varepsilon}: \mathbb{R}^{d+1} \rightarrow \mathbb{R}^{d+1} , \quad G_{y, \varepsilon}(x)=\frac{x-y}{K_\varepsilon(|x-y|)^\alpha} \, ,
    \]
    for which we have $G_{y,\varepsilon} \in C^1(\mathbb{R}^{d+1}, \mathbb{R}^{d+1})$, $|G_{y,\varepsilon}| \leq |G_y|$ and $G_{y, \varepsilon} \overset{\varepsilon \to 0}{\rightarrow} G_y$ in $C^1_{\textnormal{loc}}$ on $\mathbb{R}^{d+1} \setminus \{y\}$.
    Again we compute
    \[
    |\textnormal{div } G_{y, \varepsilon}(x)| = 
    \left| \frac{-\alpha +d+1}{K_\varepsilon (|x-y|)^{\alpha}}
    \, + \, 
    \frac{\alpha \varepsilon ^2}{K_\varepsilon (|x-y|)^{\alpha +2}} \right|
    \leq \frac{d+1}{|x-y|^{\alpha}} .
    \]
    Hence, we may again apply the Gauss-Green formula to $G_{y, \varepsilon}$ (times a suitable smooth cutoff function with compact support) and get by dominated convergence
    \begin{align*}
        &  (-\alpha+d+2) \, \int_{\Omega } 
    \int_{\partial ^\ast \Omega} 
         \frac{\langle x-y, n(x)\rangle}{|x-y|^\alpha}
        \, d \mathcal{H}^{d}(x) \, dy
        \\
        =&
        (-\alpha+d+2) \,
        \int_{\Omega } \int_{\partial ^\ast \Omega} 
        \langle G_{y}, n(x)\rangle
        \, d \mathcal{H}^{d}(x) \, dy
        \\
        =&
        (-\alpha+d+2) \,
       \lim\limits_{\varepsilon \to 0} 
        \int_{\Omega } \int_{\partial ^\ast \Omega} 
       \langle G_{y,\varepsilon}, n(x)\rangle
       \, d \mathcal{H}^{d}(x) \, dy
       \\
       =&
       (-\alpha+d+2) \,
       \lim\limits_{\varepsilon \to 0} 
       \int_{\Omega } \int_{\Omega} 
       \textnormal{div } G_{y,\varepsilon}
       \, dx \, dy
       \\
       =&
       (-\alpha+d+2) \,
       \int_{\Omega } \int_{\Omega} 
       \textnormal{div } G_{y}
       \, dx \, dy
       \\
       =&
       (-\alpha+d+2) \,
       (-\alpha +d+1) \, 
        \int_{\Omega } \int_{\Omega} 
       \frac{1}{|x-y|^{\alpha}} 
        \, dx \, dy .
    \end{align*}
    So for $\alpha <d+1$ \eqref{eq: Riesz =< TP sym} follows from \eqref{eq: Riesz TP_sym intermediate}.
    Finally, by \eqref{eq: R_m definition}, \eqref{eq: Riesz TP_sym intermediate} and dominated convergence justified by \eqref{eq: Riesz TP_sym intermediate bound} and \eqref{eq: Riesz TP_sym intermediate bound 2}
    \begin{align*}
    R_{d+1}(\Omega)=& \lim\limits_{\alpha \nearrow d+1} (-\alpha+d+1) \, R_\alpha(\Omega)
    \\
    =&
    \lim\limits_{\alpha \nearrow d+1} 
    \int_{\partial ^\ast \Omega} \int_{\Omega}
    \frac{\langle x-y, n(x)\rangle}{|x-y|^\alpha} 
    \, dy \, d \mathcal{H}^{d}(x)
    \\
    =&
    \int_{\partial ^\ast \Omega} \int_{\Omega}
    \frac{\langle x-y, n(x)\rangle}{|x-y|^{d+1}} 
    \, dy \, d \mathcal{H}^{d}(x)
    \\
    =&
    \int_{\partial ^\ast \Omega}
    \int_{\partial ^\ast \Omega}
    \frac{\langle x-y, n(x)\rangle \, \langle y-x, n(y) \rangle}{|y-x|^{d+1}} 
    \, d \mathcal{H}^{d}(y) \, d \mathcal{H}^{d}(x)  .
    \end{align*}
\end{proof}

\begin{rmk}
\label{rmk: TP_sym and fractional perimeter relation}
An argument similar to that in \Cref{lem: TP_sym and Riesz relation} yields for $\alpha \in (d+1, d+2)$
\begin{align*}
&\int_{\Omega} \int_{\Omega^C}
\frac{1}{|x-y|^{\alpha}} \, dx \, dy 
\\
=&
 \frac{-1}{(-\alpha+d+2)(-\alpha+d+1)}
	 \, 
        \int_{\partial ^\ast \Omega}
        \int_{\partial ^\ast \Omega}
        \frac{\langle x-y, n(x)\rangle \, \langle y-x, n(y) \rangle}{|y-x|^{\alpha}} 
        \, d \mathcal{H}^{d}(y) \, d \mathcal{H}^{d}(x)  .
\end{align*}
 Here the left-hand side corresponds to the common definition of the fractional perimeter, which is known to be uniquely minimized by balls upon volume fixation \cite[Equation~(4.2)]{FracPerIso}. However, as the sharp upper bound on the Riesz energy (see \Cref{prop: Riesz maximum}) is crucial for our argumentation, we will not exploit this relation any further. 
\end{rmk}

The sharp inequality \eqref{eq: Riesz maximum} for $R_\alpha$ further translates to the chord power integral $I_\alpha$ defined in \eqref{eq: def I_alpha} via the following lemma. 
\begin{lem}
\label{lem: sharp lower bound on Riesz energy}
    Let $d \geq 1$, $\alpha \in [0,d+1]$ and $\Omega \subset \mathbb{R}^{d+1}$ be a bounded Caccioppoli set. For $x \in \partial ^\ast \Omega$ let $S_x:=\{\omega \in \mathbb{S}^{d} \, | \, \langle \omega , n_\Omega(x) \rangle <0 \}$ and $l(x, \omega)$ as in \Cref{nota: L and l}. Then for the Riesz energy \eqref{eq: Riesz energy} there holds
    \begin{equation}
    \label{eq: sharp lower bound on Riesz energy}
     c_{d, \alpha} \, R_\alpha (\Omega) \geq 
    \int_{\partial ^\ast \Omega} \int_{S_x} 
          l(x, \omega) ^{-\alpha+d+2} \, \langle -\omega , n_\Omega(x) \rangle
         \, d \mathcal{H}^{d}(\omega) \, d \mathcal{H}^{d}(x) \, 
    \end{equation}
    with $c_{d, \alpha}= (-\alpha +d+2)(- \alpha+d+1)$ for $\alpha <d+1$ and $c_{d,d+1}=1$. If $\alpha <d+1$, equality holds if and only if $\partial ^\ast \Omega$ is the reduced boundary of a convex set. If $\alpha=d+1$, then \eqref{eq: sharp lower bound on Riesz energy} is an equality.
\end{lem}

\begin{proof}
    Recall $L_{x, \omega}$ and $l(x, \omega)$ from \Cref{nota: L and l}.
    By \Cref{lem: finite perimeter sets satisfy SOCC} for $\mathcal{H}^{d}$-a.e. $x \in \partial ^\ast \Omega$ and $\mathcal{H}^{d}$-a.e. $\omega \in \mathbb{S}^{d}$ there exists a finite family $t_1<t_2<...<t_{2k}$ such that $L_{x, \omega} \cap \partial ^\ast \Omega =\{x+ t_i \omega \, | \, 1 \leq i \leq 2k\}$ and 
    \[
    \textnormal{sgn } \langle \omega , n (x+ t_i \omega) \rangle = (-1)^i  .
    \]
    Fix such an $x \in \partial ^\ast \Omega$. Let $i^\ast \in \{1,..., 2k\}$ denote the index such that $t_{i^\ast}=0$.
    By \eqref{eq: TP_ssym inner integral finite} we may apply \Cref{prop: coarea formula} to $U=\partial ^\ast \Omega \setminus \{x\}$, $V= \mathbb{S}^{d}$, $f=\phi_x$ and 
    \[
    g: \partial ^\ast \Omega \setminus \{x\} \rightarrow \R , \quad 
    g(y):= \langle -\phi_x(y), n(x) \rangle  \, |x-y|^{-\alpha+d+2} \, \textnormal{sgn} \langle \phi_x(y), n(y) \rangle   .
    \]
    Then we may write
    \begin{align*}
        &\int_{\partial ^\ast \Omega}
        \frac{\langle x-y, n(x)\rangle \, \langle y-x, n(y) \rangle}{|y-x|^\alpha} 
        \, d  \mathcal{H}^{d}(y)
        \\
        =&
        \int_{\partial ^\ast \Omega}
         \langle -\phi_x(y), n(x) \rangle  \, |x-y|^{-\alpha+d+2} \, \textnormal{sgn} \langle \phi_x(y), n(y) \rangle  \, |J\phi_x(y)| 
         \, d  \mathcal{H}^{d}(y)
        \\
        =&
        \int_{\mathbb{S}^{d}} 
       \sum_{\phi_x^{-1} (\omega)} 
       \langle - \omega , n(x) \rangle  \, |x-y|^{- \alpha +d+2} \,
       \textnormal{sgn} \langle \omega, n(y) \rangle 
         \, d  \mathcal{H}^{d}(\omega)
       \\
       =&
         \int_{\mathbb{S}^{d}} 
         \langle - \omega , n(x) \rangle
      \sum_{\overset{{\lambda >0}}{x+\lambda \omega \in \partial ^\ast \Omega}}
         \, \lambda^{- \alpha +d+2} \,
       \textnormal{sgn} \langle \omega, n(x+\lambda \, \omega) \rangle
        \, d  \mathcal{H}^{d}(\omega)
       \\
        =&
      \int_{\mathbb{S}^{d}}
        \langle -\omega , n(x) \rangle 
      \sum_{i=i^\ast+1}^{2k}
         \, t_i^{- \alpha +d+2} \,
       \textnormal{sgn} \langle \omega, n(x+t_i \, \omega) \rangle
        \, d  \mathcal{H}^{d}(\omega) 
       \\
       =&
       \int_{\mathbb{S}^{d}}
       \langle -\omega , n(x) \rangle 
       \sum_{i=i^\ast+1}^{2k}
       \, t_i^{- \alpha +d+2} \,
       (-1)^{i}
         \, d  \mathcal{H}^{d}(\omega)
       \\
       =&
      \int_{\mathbb{S}^{d}}
       |\langle \omega , n(x) \rangle| \, (-1)^{i^\ast +1} 
       \sum_{i=i^\ast+1}^{2k}
       \, t_i^{- \alpha +d+2} \,
       (-1)^{i}
        \, d  \mathcal{H}^{d}(\omega) .
    \end{align*}
    Integrating with respect to $x$ yields together with \Cref{lem: TP_sym and Riesz relation} the identity
    \[
    c_{d,\alpha} R_\alpha(\Omega)= 
    \int_{\partial ^\ast \Omega} \int_{\mathbb{S}^{d}}
    |\langle \omega , n(x) \rangle| \, (-1)^{i^\ast +1} 
    \sum_{i=i^\ast+1}^{2k}
    \, t_i^{- \alpha +d+2} \,
    (-1)^{i}
    \, d \mathcal{H}^{d}(\omega) \, d \mathcal{H}^{d}(x) \, ,
    \]
    where the right-hand side exists a priori only as an iterated integral.
    However, note that further by \Cref{lem: finite perimeter sets satisfy SOCC}, for $\mathcal{H}^{d}$-a.e. $x \in \partial ^\ast \Omega$ and $\mathcal{H}^{d}$-a.e. $\omega \in \mathbb{S}^{d}$ there holds $\chi_{\Omega}(x+t\omega)=1$ if and only if $t\in [t_i,t_{i+1}]$ and $\langle \omega , n(x+t_i \omega) \rangle <0$ (up to an $\mathcal{L}^1$-null set). Therefore,
    \[
    \sum_{i=i^\ast+1}^{2k}
    \, t_i^{- \alpha +d+2} \,
    (-1)^i
    =
    (- \alpha + d+2) \, \int_0^{\infty} t^{-\alpha + d+1} \chi_{\Omega }(x+t \omega) \, dt  \leq \textnormal{diam}(\Omega)^{d+2-\alpha} \, ,
    \]
    so the above integral exists as an integral on the product space $\Sigma ^\ast \times \mathbb{S}^{d}$.
    
    Fix $x, \omega$ satisfying the above property, and denote by 
    \[
    \pi_\omega: \Sigma ^\ast  \rightarrow V_\omega , \quad \pi_\omega (y):=y- \langle y, \omega \rangle \omega 
    \]
    the orthogonal projection onto $V_\omega:=\omega^\perp \cap B_R(0)$. Here, $R>0$ is chosen sufficiently large that the map is well defined. Then for $z=\pi_\omega(x)$ and $s_i=t_i+\langle x, \omega \rangle $ there holds 
    \[
    L_{x,\omega} \cap \partial ^\ast \Omega=L_{z,\omega} \cap \partial ^\ast \Omega=\{ z+ s_i \omega \, | \, 1 \leq i \leq 2k \} \, , \quad  x=z+s_{i^\ast} \omega  \quad \textnormal{and} \quad t_i=s_i-s_{i^\ast} \, .
    \] 
    Setting $J:= \bigcup _{j=1}^{k} (s_{2j-1}, s_{2j})$, we also have
    \[
    \chi_{\Omega}(z+s\omega)= \chi_J(s) \, \textnormal{for a.e. } s \in \mathbb{R} \, 
    \]
    and 
    \begin{equation}
    	\label{eq: l via s}
    l(x, \omega) = \min\limits_{y \in \phi_x^{-1}(\omega)} |x-y|
    =
    (s_{i^\ast +1}-s_{i^\ast}) \quad \textnormal{if } \omega \in S_x \, ,
    \end{equation}
    while $l(x, \omega)=0$ if $\omega \in \mathbb{S}^{d} \setminus \bar{S}_x$.
    Thus, via Fubini's theorem and \Cref{prop: coarea formula} applied to $U=\partial ^\ast \Omega$, $V=V_\omega$, $f=\pi_ \omega$ with $ |J \pi_\omega (x) | = |\langle \omega , n(x) \rangle|$ and 
    \[
    g=g_\omega: \partial ^\ast \Omega \rightarrow \R, \, \quad 
    g_\omega(x):= (-1)^{i^\ast +1} 
    \sum_{i=i^\ast+1}^{2k}
    \, t_i^{- \alpha +d+2} \,
    (-1)^{i} \, ,
    \]
    we get
    \begin{align*}
     &\int_{\partial ^\ast \Omega}     \int_{\mathbb{S}^{d}}
    |\langle \omega , n(x) \rangle| \, (-1)^{i^\ast +1} 
    \sum_{i=i^\ast+1}^{2k}
    \, t_i^{- \alpha +d+2} \,
    (-1)^{i}
     \, d  \mathcal{H}^{d}(\omega) \, d \mathcal{H}^{d} (x)
    \\
    =&
    \int_{\mathbb{S}^{d}}  \int_{\partial ^\ast \Omega}
    |\langle \omega , n(x) \rangle| 
    \sum_{i=i^\ast+1}^{2k}
    \, t_i^{- \alpha +d+2} \, 
    (-1)^{i+i^\ast +1}
     \, d  \mathcal{H}^{d}(x) \, d \mathcal{H}^{d} (\omega)
    \\
    =&
    \int_{\mathbb{S}^{d}} \int_{V_\omega}
    \sum_{\pi_\omega^{-1}(z)}
    \sum_{i=i^\ast+1}^{2k}
    \, t_i^{- \alpha +d+2} \,
    (-1)^{i+i^\ast +1}
    \, d  \mathcal{H}^{d}(z) \, d \mathcal{H}^{d} (\omega)
    \\
    =&
    \int_{\mathbb{S}^{d}} \int_{V_\omega}
    \sum_{x \in L_{z, \omega} \cap \partial ^\ast \Omega}
    \sum_{i=i^\ast+1}^{2k}
    (s_i-s_{i^\ast})^{-\alpha +d+2} \, 
    (-1)^{i+i^\ast +1}
    \, d  \mathcal{H}^{d}(z) \, d \mathcal{H}^{d} (\omega)
    \\
    =&
    \int_{\mathbb{S}^{d}} \int_{V_\omega}
    \sum_{1 \leq j <i \leq 2k}
    (s_i-s_j)^{-\alpha + d+2} \, (-1)^{i+j+1}
   \, d  \mathcal{H}^{d}(z) \, d \mathcal{H}^{d} (\omega) ,
    \end{align*}
    where we have set $k=0$ if $\pi_\omega ^{-1}(z) =\emptyset$. For fixed $\omega \in \mathbb{S}^{d}$ and $z \in V_\omega$ (if possible) let $s_1<s_2<...<s_{2k}$, $k \in \mathbb{N}_0$ be given as in the previous discussion. 
    If $\alpha <d+1$ we conclude from \Cref{lem: convexity inequality} by strict convexity of $s \mapsto |s|^{-\alpha+d+2}$ that 
    \[
    \sum_{1 \leq j <i \leq 2k}
      (s_i-s_j)^{-\alpha + d+2} \, (-1)^{i+j+1} 
      \geq 
      \sum_{j=1}^{k} (s_{2j}-s_{2j-1})^{-\alpha + d+2}  
    \]
    with equality if and only if $k \in \{0,1\}$. Otherwise, if $\alpha=d+1$ then $s \mapsto s^{-\alpha +d+2}$ is the identity map, so by \Cref{lem: convexity inequality} 
    \[
    \sum_{1 \leq j <i \leq 2k}
    (s_i-s_j)^{-\alpha + d+2} \, (-1)^{i+j+1} 
    = 
    \sum_{j=1}^{k} (s_{2j}-s_{2j-1})^{-\alpha + d+2}  .
    \]
    Thus, defining $A_\omega(z):= \{x \in \pi_\omega ^{-1}(z) \, | \, \omega \in S_x \}$ we obtain
    \begin{align*}
        & \int_{\mathbb{S}^{d}} \int_{V_\omega}
      \sum_{1 \leq j <i \leq 2k}
      (s_i-s_j)^{-\alpha + d+2} \, (-1)^{i+j+1}
      \, d  \mathcal{H}^{d}(z) \, d \mathcal{H}^{d} (\omega)
        \\
        \geq&
        \int_{\mathbb{S}^{d}} \int_{V_\omega}
        \sum_{j=1}^{k} (s_{2j}-s_{2j-1})^{-\alpha + d+2}
       \, d  \mathcal{H}^{d}(z) \, d \mathcal{H}^{d} (\omega)
        \\
        =&
        \int_{\mathbb{S}^{d}} \int_{V_\omega}
        \sum_{x \in A_\omega (z) } 
        (s_{i^\ast +1} - s_{i^\ast})^{-\alpha + d+2}
       \, d  \mathcal{H}^{d}(z) \, d \mathcal{H}^{d} (\omega) 
    \end{align*}
    with equality if and only if $\alpha =d+1$ or for $\mathcal{H}^{d}$-a.e. $\omega \in \mathbb{S}^{d}$ and $\mathcal{H}^{d}$-a.e. $z \in \omega ^\perp$ there holds $k\in \{0,1\}$. As the succeeding argument shows, the latter equality condition implies $k=1$ for $\mathcal{H}^{d}$-a.e. $\omega \in \mathbb{S}^{d}$ and $\mathcal{H}^{d}$-a.e. $x \in \partial ^\ast \Omega$. By \cite[Remark~15.2]{FinitePerimeterSets} and the proof of \cite[Corollary, p.~4]{bushling2024singular} this further implies that $\partial ^\ast \Omega$ is the reduced boundary of a convex set.
    
    For $\omega \in \mathbb{S}^{d}$ define $Q_\omega:= \{ x \in \partial ^\ast \Omega \, | \, k>1\}$ and further $Q=\{(\omega, x) \in \mathbb{S}^{d} \times \partial ^\ast \Omega \, | \, x \in Q_\omega \}$. Suppose $\mathcal{H}^{2d}(Q)>0$. Then by \Cref{prop: coarea formula} and $\mathcal{H}^{2d}$-a.e. positivity of $|J \pi_\omega (x) |=|\langle \omega, n(x) \rangle |$ (see \cite[Theorem, p.~9]{bushling2024singular}) we have
    \[
    0< \int_Q |J \pi_\omega (x)| 
    \, d \mathcal{H}^{d}(x) \, d \mathcal{H}^{d}(\omega)
    =
    \int_{\mathbb{S}^{d}} \int_{V_\omega} \mathcal{H}^{0}(\pi_\omega ^{-1} (z) \cap Q_\omega)
    \, d \mathcal{H}^{d}(z) \, d \mathcal{H}^{d}(\omega) .
    \]
    Hence, also $\mathcal{H}^{2d}(\{ (\omega, z) \in \mathbb{S}^{d} \times \omega ^\perp \, | \, k>1 \})>0$.
    
    Using \eqref{eq: l via s} and applying the coarea formula and Fubini's theorem once again, we finally get 
    \begin{align*}
    	c_{d,\alpha} \, R_{\alpha}(\Omega) \geq
        &\int_{\mathbb{S}^{d}} \int_{V_\omega}
        \sum_{x \in A_\omega (z)} 
        (s_{i^\ast +1} - s_{i^\ast})^{-\alpha + d+2} 
       \, d  \mathcal{H}^{d}(z) \, d \mathcal{H}^{d} (\omega)
        \\
        =&
        \int_{\mathbb{S}^{d}} \int_{V_\omega}
        \sum_{x \in \pi_\omega ^{-1}(z)} 
        l(x, \omega)^{-\alpha + d+2} 
        \, d  \mathcal{H}^{d}(z) \, d \mathcal{H}^{d} (\omega)
        \\
        =&
       \int_{\mathbb{S}^{d}} \int_{\partial ^\ast \Omega} 
       l(x, \omega)^{-\alpha + d+2} \, |\langle \omega , n(x) \rangle|
      \, d  \mathcal{H}^{d}(x) \, d \mathcal{H}^{d} (\omega)
       \\
       =&
       \int_{\partial ^\ast \Omega} \int_{S_x} 
       l(x, \omega) ^{-\alpha+d+2} \, \langle -\omega , n(x) \rangle
      \, d  \mathcal{H}^{d}(\omega) \, d \mathcal{H}^{d} (x) 
    \end{align*}
    with equality if and only if $\alpha=d+1$ or $\partial ^\ast \Omega$ is the reduced boundary of a convex set.
\end{proof}

We are now prepared to prove our second main result.

\begin{proof}[Proof of \Cref{thm: main 2}]
Let $d,q,p , \Sigma, \Omega$ be given as above. If $p=q+d$, then \eqref{eq: TP^(p,q) minimum sets 1} follows directly from \Cref{thm: lower bound on TP_q+d for sets}, subsequently implying \eqref{eq: TP^(p,q) minimum sets 2} via the classical isoperimetric inequality (see \Cref{rmk: implications of quotient bounds}). Therefore, in the following we assume $p>q+d$. Moreover, w.l.o.g. we assume $\TP^{(p,q)}(\Sigma)< \infty$. In order to show \eqref{eq: TP^(p,q) minimum sets 1} and its equality cases, we may assume that $S=S_1=\partial B_1$ and $\mathcal{H}^d(\Sigma ^\ast)=\mathcal{H}^{d}(S_1)=1$; for \eqref{eq: TP^(p,q) minimum sets 2} and its equality cases, we may assume that $S=S_2=\partial B_2$ and $\mathcal{L}^{d+1}(\Omega)=\mathcal{L}^{d+1}(B_2)=1$. It therefore suffices to show under the respective assumption that $\TP^{(p,q)}(\Sigma) \geq \TP^{(p,q)}(S_i)$ with equality if and only if $\Sigma ^\ast= S_i$ (up to translations).
   
As in \Cref{thm: lower bound on TP_q+d for sets}, we first bound $\TP^{(p,q)}(\Sigma)$ from below using Young's inequality in the following way
    \begin{align*}
        &\TP^{(p,q)}(\Sigma) 
        \\
        =& 
        \frac{q-1}{q} 
        \int_{\Sigma ^\ast} \int_{\Sigma ^\ast}
        \frac{|\langle x-y, n(y) \rangle |^q}{|x-y|^p} 
        \, d \mathcal{H}^{d}(y) \,  d \mathcal{H}^{d}(x)
        \\
        +&
        \frac{1}{q}
        \int_{\Sigma ^\ast} \int_{\Sigma ^\ast}
        \frac{|\langle x-y, n(x) \rangle |^q}{|x-y|^p} 
        \, d \mathcal{H}^{d}(y) \,  d \mathcal{H}^{d}(x)
        \\
        =&
        \int_{\Sigma ^\ast} \int_{\Sigma ^\ast}
        \frac{\frac{q-1}{q}| \langle x-y, n(x) \rangle |^q \, + \, \frac{1}{q} | \langle x-y, n(y) \rangle |^q}{|x-y|^p} 
        \, d \mathcal{H}^{d}(y) \,  d \mathcal{H}^{d}(x)
        \\
       \geq&
       \int_{\Sigma ^\ast} \int_{\Sigma ^\ast}
       \frac{|\langle x-y, n(x) \rangle |^{q-1} \, 
       |\langle x-y, n(y) \rangle |}{|x-y|^p} 
       \, d \mathcal{H}^{d}(y) \,  d \mathcal{H}^{d}(x) .
    \end{align*}
    Here, equality holds for spheres $S$, since 
    \[
    |\langle x-y, n(x) \rangle|= |\langle x-y, n(y) \rangle| \quad \forall x,y \in S  .
    \]
    Applying \Cref{prop: coarea formula} to $U=\Sigma ^\ast \setminus \{x\}$, $V=\mathbb{S}^{d}$, $g=\frac{|\langle x- \cdot,n(x) \rangle|^{q-1}}{|x-\cdot|^{p-d-1}}$ and $f=\phi_x$ we further get
    \begin{align*}
    	\TP^{(p,q)}(\Sigma)\geq
        &\int_{\Sigma ^\ast} \int_{\Sigma ^\ast}
       \frac{|\langle x-y, n(x) \rangle |^{q-1} \, 
       |\langle x-y, n(y) \rangle |}{|x-y|^p} 
       \, d \mathcal{H}^{d}(y) \,  d \mathcal{H}^{d}(x)
       \\
       =&
       \int_{\Sigma ^\ast} \int_{\Sigma ^\ast}
       |\langle \phi_x(y), n(x) \rangle |^{q-1} \, |x-y|^{q+d-p} \, |J \phi_x(y) |
       \, d \mathcal{H}^{d}(y) \,  d \mathcal{H}^{d}(x)
       \\
       =&
       \int_{\Sigma ^\ast} \int_{\mathbb{S}^{d}} 
       \sum_{\phi_x^{-1} (\omega)} 
       |\langle \omega , n(x) \rangle |^{q-1} \, |x-y|^{q+d-p} 
       \, d \mathcal{H}^{d}(\omega) \,  d \mathcal{H}^{d}(x)
       \\
       \geq&
        \int_{\Sigma ^\ast} \int_{S_x} 
       \langle - \omega , n(x) \rangle ^{q-1}
       \sum_{\phi_x^{-1} (\omega)} 
        |x-y|^{q+d-p} 
       \, d \mathcal{H}^{d}(\omega) \,  d \mathcal{H}^{d}(x)
        \\
       \geq&
       \int_{\Sigma ^\ast} \int_{S_x} 
       \langle - \omega , n(x) \rangle ^{q-1}
       (\min_{\phi_x^{-1} (\omega)} 
       |x-y|)^{q+d-p} 
       \, d \mathcal{H}^{d}(\omega) \,  d \mathcal{H}^{d}(x) \, ,
    \end{align*}
    with equality for spheres $S$, since then $\phi_x$ is injective for all $x \in S$ and $\textnormal{Im}(\phi_x)=S_x$.
    Recalling $l(x,\omega)$ from \Cref{nota: L and l}, by \eqref{eq: l via s} we may rewrite the lower bound as
    \begin{align*}
    	 &\int_{\Sigma ^\ast} \int_{S_x} 
    	\langle - \omega , n(x) \rangle ^{q-1}
    	(\min_{\phi_x^{-1} (\omega)} 
    	|x-y|)^{q+d-p} 
    	\, d \mathcal{H}^{d}(\omega) \,  d \mathcal{H}^{d}(x)
    	\\
    	=&
	    \int_{\Sigma ^\ast} \int_{S_x} 
	    \langle - \omega , n(x) \rangle ^{q-1} \, l(x,\omega)^{q+d-p}
	    \, d \mathcal{H}^{d}(\omega) \,  d \mathcal{H}^{d}(x)  .
	\end{align*}
    Notice that for spheres $S$ we have $l(x, \phi_x(y))=|x-y|$ for all $x,y \in S, x\neq y$. Denote by $D_i$ the diameter of $S_i$, $i \in \{1,2\}$. Then by \Cref{prop: Supporting Hyperplane Theorem strong} and \Cref{cor: sphere constant integrand} there holds
    \[
    \langle -\phi_x (y), n(x) \rangle=
    |\langle \frac{x-y}{|x-y|}, n(x) \rangle|
    =
    D_i^{-1} \, |x-y| \quad \forall x,y \in S_i, \, x \neq y  .
    \]
    Equivalently,
    \[
    D_i \, \langle -\omega, n(x) \rangle =
    l(x, \omega) \quad \forall x \in S_i, \omega \in S_x \, .
    \]
    Setting $\kappa:=2q+d-2-p>0$, we therefore have by \Cref{lem: Young's inequality} applied to $ a=(D_i \,  \langle -\omega, n(x) \rangle)^{q-2}$, $b=l(x,\omega)^{q+d-p}$,  $\varrho_1= \frac{\kappa}{q-2}$ and $\varrho_2=-\frac{\kappa}{q+d-p}$ 
    \begin{align*}
    	D_i^{q-2} \TP^{(p,q)}(\Sigma) \geq
        &\int_{\Sigma ^\ast} \int_{S_x} 
       (D_i \,\langle - \omega , n(x) \rangle) ^{q-2} \,  l(x,\omega)^{q+d-p}
       \, \langle - \omega , n(x) \rangle
       \, d \mathcal{H}^{d}(\omega) \,  d \mathcal{H}^{d}(x)
       \\
      \geq&
       \frac{q-2}{\kappa } D_i^\kappa
       \int_{\Sigma ^\ast} \int_{S_x} 
       \langle - \omega , n(x) \rangle^{\kappa +1}
       \, d \mathcal{H}^{d}(\omega) \,  d \mathcal{H}^{d}(x)
       \\
       +&
       \frac{q+d-p}{\kappa }
       \int_{\Sigma ^\ast} \int_{S_x} 
       l(x,\omega)^{\kappa} \, \langle -\omega , n(x) \rangle
       \, d \mathcal{H}^{d}(\omega) \,  d \mathcal{H}^{d}(x) \, ,
    \end{align*}
    where equality holds for $\Sigma ^\ast= S_i$.
   The first inner integral is invariant under rotations, so we may replace $n(x)$ by $e_1$ and $S_x$ by $A:=\{ v \in \mathbb{S}^{d} \, | \, \langle v, e_1 \rangle <0 \}$, yielding
   \begin{align*}
 	\int_{S_x} 
 	\langle - \omega , n(x) \rangle^{\kappa +1}
 	\, d \mathcal{H}^{d}(\omega)
    =
   	\int_{A} 
   	\langle - \omega , e_1 \rangle^{\kappa +1}
   	\, d \mathcal{H}^{d}(\omega)
   	= 
   	\int_{\mathbb{S}^{d}} 
   	\frac{|\langle  \omega , e_1 \rangle|^{\kappa +1}}{2}
   	\, d \mathcal{H}^{d}(\omega) = b_{d,\kappa+2}
   \end{align*}
   with $b_{d,\kappa+2}$ from \Cref{cor: lower bound on F_q for sets}. Moreover, since for $\alpha:= d+2- \kappa$ we have $\alpha \in [0,d+1]$, according to \Cref{lem: sharp lower bound on Riesz energy} the integral in the second summand is a sharp lower bound of $c_{d, \alpha} R_{\alpha}(\Omega)$. Thus, by negativity of $\frac{q+d-p}{\kappa}$
   \[
   	D_i^{q-2} \TP^{(p,q)}(\Sigma) \geq
   \frac{q-2}{\kappa} \, D_i^\kappa \, b_{d,\kappa +2} \mathcal{H}^{d}(\Sigma ^\ast)
   \, + \,
   \frac{q+d-p}{\kappa}
   c_{d, \alpha} R_\alpha(\Omega)  .
   \] 
   \Cref{prop: Riesz maximum} applied to the second summand then yields
    \[
    	D_i^{q-2} \TP^{(p,q)}(\Sigma) 
      \geq
       \frac{q-2}{\kappa} \, D_i^\kappa \, b_{d,\kappa +2} \mathcal{H}^{d}(\Sigma ^\ast)
       \, + \,
       \frac{q+d-p}{\kappa}
       c_{d, \alpha} R_\alpha(B_i) \frac{\mathcal{L}^{d+1}(\Omega)^{\frac{2d+2-\alpha}{d+1}}}{\mathcal{L}^{d+1}(B_i)^{\frac{2d+2-\alpha}{d+1}}}  .
    \]
    Equality holds in both inequalities above, if $\Sigma^\ast=S_i$. 
    
    Finally, by the classical isoperimetric inequality (see \cite[Theorem~14.1]{FinitePerimeterSets})
    \[
    	D_i^{q-2} \TP^{(p,q)}(\Sigma) \geq
     \frac{q-2}{\kappa} \, D_i^\kappa \, b_{d,\kappa +2} \mathcal{H}^{d}(S_i)
     \, + \,
     \frac{q+d-p}{\kappa}
     c_{d, \alpha} R_\alpha(B_i) \, ,
    \]
    where equality holds for $\Sigma^\ast=S_i$, so in particular the right-hand side is equal to $D_i^{q-2} \TP^{(p,q)}(S_i)$. On the other hand, if equality holds above, then by the classical isoperimetric inequality there exist $z \in \mathbb{R}^{d+1}$ and $r>0$ such that the symmetric difference between $\Omega$ and $B_r(z)$ has Lebesgue measure zero, which by \cite[Remark~15.2]{FinitePerimeterSets} implies that $\Sigma ^\ast=S_i$ (up to translations). 
\end{proof}

\begin{rmk}
\label{rmk: implications of quotient bounds}
    Depending on the sign of $h:=q+2d-p$, as a consequence of the classic isoperimetric inequality \cite[Theorem~14.1]{FinitePerimeterSets}, a sharp bound on either one of the two quotients in \eqref{eq: scale invariant quotients} via the value at round hyperspheres implies the other. For example, if $h<0$, then the latter bound implies the former, since
    \[
    \frac{\TP^{(p,q)}(\Sigma)}{\mathcal{H}^d(\Sigma ^\ast)^{\frac{h}{d}}} 
    \geq
    \omega _d \, \frac{\TP^{(p,q)}(\Sigma)}{ \mathcal{L}^{d+1}(\Omega)^{\frac{h}{d+1}}} 
    \geq 
    \omega _d \, \frac{\TP^{(p,q)}(S)}{\mathcal{L}^{d+1}(B)^{\frac{h}{d+1}}} 
    =
    \frac{\TP^{(p,q)}(S)}{\mathcal{H}^{d}(S)^{\frac{h}{d}}} \, ,
    \]
    where $\omega _d := \mathcal{L}^{d+1}(B)^{\frac{h}{d+1}} / \mathcal{H}^{d}(S)^{\frac{h}{d}}$.
\end{rmk}

\begin{rmk}
    \Cref{thm: main 2} extends the parameter range for which \cite[Theorem~1.1]{döhrer2026fencheltheoremgaussmaps} holds for planar simple closed curves ($d=1$) from 
    \[
    q \geq 1, p \in [q+1, 2q+1) \cap [2q-2,4q-2]
    \]
    to 
    \[
    q \geq 1, p \in [q+1, 2q+1) \cap [2q-4,4q-2] .
    \]
\end{rmk}

\section*{Acknowledgments}
ChatGPT helped us to find the convex analysis argument of \Cref{lem: convexity inequality} to prove \Cref{lem: sharp lower bound on Riesz energy}. Other than that, AI-tools were only used for proof checking of the manuscript.
The authors take full responsibility for the mathematical correctness and final content of the paper.
The authors would also like thank Wenhui Shi for her helpful comments, which ultimately lead to an extension of the domain in \Cref{thm: main 2} from boundaries of convex sets to boundaries of arbitrary finite-perimeter sets.
AD is funded by the DFG-Graduiertenkolleg Energy, Entropy, and Dissipative Dynamics (EDDy), project no. 320021702/GRK2326. HvdM gratefully acknowledges the support by the Excellence Initiative of the German federal and state governments.

\appendix
\section*{Appendix}
\addcontentsline{toc}{section}{Appendix}

\setcounter{thm}{0}
\renewcommand{\thethm}{A.\arabic{thm}}

\begin{lem}
	\label{lem: convexity inequality}
	Let $f:\mathbb{R} \rightarrow \mathbb{R}$ be strictly convex or the identity map id, $k \in \mathbb{N}_0$ and $x_1<x_2<...<x_{2k}$. Then
	\[
	\sum_{1\leq j <i \leq 2k} f(x_i-x_j) \, (-1)^{i+j+1} \geq \sum_{j=1}^k f(x_{2j}-x_{2j-1}) 
	\]
	with equality if and only if $k \in \{0,1\}$ or $f=\textnormal{id}$.
\end{lem}

\begin{proof}
	We have
	\begin{align*}
		&\sum_{1\leq j <i \leq 2k} f(x_i-x_j) \, (-1)^{i+j+1} 
		\, - \, 
		\sum_{j=1}^k f(x_{2j}-x_{2j-1})
		\\
		=&
		\sum_{1\leq j <i \leq k}
		\Bigg[
		\Big( f(x_{2i}-x_{2j-1}) \, - \, f(x_{2i-1}-x_{2j-1}) \Big)
		\, - \, 
		\Big( f(x_{2i}-x_{2j}) \, - \, f(x_{2i-1} -x_{2j}) \Big)
		\Bigg]  .
	\end{align*}
	If $k \in \{0,1 \}$, then the above sum is empty, and if $f=\textnormal{id}$, then each summand is equal to zero, so it remains to consider $k >1$ and strictly convex $f$.
	It suffices to show that each summand is positive.
	By strict convexity of $f$, finite differences are strictly monotonically increasing, so in particular for $a < b$ and $h > 0$ there holds
	\[
	f(a+h)-f(a) < f(b+h)-f(b)  .
	\]
	Applying the above inequality to 
	\[
	a=x_{2i-1} - x_{2j} \,  , \quad b=x_{2i-1}- x_{2j-1} \, , \quad h=x_{2i}- x_{2i-1}
	\]
	completes the proof.
\end{proof}

\begin{lem}[Variant of Young's inequality] 
	\label{lem: Young's inequality}
	Let $p \in (0,1)$, $q \in (0, \infty)$ satisfy $\frac{1}{p}-\frac{1}{q}=1$. Then for all $a \geq 0, b>0$ one has
	\[
	ab \geq \frac{a^p}{p} \, - \, \frac{b^{-q}}{q}  
	\]
	with equality if $a^p=b^{-q}$.
\end{lem}
\begin{proof}
	The claim follows by applying Young's inequality $xy \leq \frac{x^u}{u} + \frac{y^v}{v}$ to $u=\frac{1}{p}$, $v=\frac{q}{p}$ (which satisfy $\frac{1}{u}+ \frac{1}{v}=1$) and $x=(ab)^p$, $y=b^{-p}$ and rearranging the terms.
\end{proof}

\newpage
\printbibliography[heading=bibintoc]
\thispagestyle{plain}

\end{document}